\documentclass[aop,noinfoline%,preprint,authoryear
]{imsart}
\usepackage{amscd,amsfonts,amssymb,amsmath,amsthm,latexsym,array,hhline,mathrsfs,booktabs,fancybox,calc,textcomp,xcolor,natbib,graphicx,parskip%,bbding
}
\usepackage[sans]{dsfont}
\usepackage[latin1]{inputenc}
\setattribute{journal}{name}{}
\numberwithin{equation}{section}

\newtheorem{theorem}{Theorem}[section]
\newtheorem{lemma}[theorem]{Lemma}
\newtheorem{proposition}[theorem]{Proposition}
\newtheorem{corollary}[theorem]{Corollary}

\newenvironment{remark}[1][Remark.]{\begin{trivlist}
\item[\hskip \labelsep {\scshape #1}]}{\end{trivlist}}

\def\E{\mathbf{E}}

\def\N{\mathbb{N}}

\def\R{\mathbb{R}}

\def\GG{\mathcal{G}}

\def\NN{\mathcal{N}}
\def\XX{\mathcal{X}}
\def\CC{\mathcal{C}}

\def\SS{\mathcal{S}}

\def\N{\mathbb{N}}
\def\P{\mathbf{P}}

\renewcommand{\d}{\mathrm{d}}
\newcommand{\e}{\mathrm{e}}
\newcommand{\id}{\mathrm{Id}}

\renewcommand{\AA}{\mathcal{A}}
\newcommand{\BB}{\mathcal{B}}
\newcommand{\DD}{\mathcal{D}}

\renewcommand{\GG}{\mathcal{G}}

\newcommand{\tr}{\mathrm{tr}}

\renewcommand{\hat}{\widehat}
\renewcommand{\Upsilon}{Z}
\renewcommand{\tilde}{\widetilde}
\begin{document}
\sloppy
\begin{frontmatter}

\title{Accuracy of empirical projections of high-dimensional Gaussian matrices}%\protect\thanksref{T1}}
\runtitle{Accuracy of empirical projections}
%\thankstext{T1}{Footnote to the title with the `thankstext' command.}
%
\begin{aug}
  \author{\fnms{Angelika}  \snm{Rohde}%\corref{}\thanksref{t2}
\ead[label=e1]{angelika.rohde@math.uni-hamburg.de}}

  \runauthor{A. Rohde}
  \affiliation{Universit\"at Hamburg}
  \address{Universit\"at Hamburg, Fachbereich Mathematik, Bundesstra{\ss}e 55, 20146 Hamburg, Germany\\
          \printead{e1}}
%
%  \address{Address of the Third author,\\
%          \printead{e3,u1}}
%
\end{aug}
\small
\begin{abstract}\scriptsize
\noindent
Let $X=C+\mathrm{E}$  with a deterministic matrix $C\in\R^{M\times M}$ and $\mathrm{E}$ some centered Gaussian $M\times M$-matrix whose entries are independent with variance $\sigma^2$.  In the present work, the accuracy of reduced-rank projections of $X$ is studied.
Non-asymptotic universal upper and lower bounds are derived, and favorable and unfavorable prototypes of matrices $C$ in terms of the accuracy of approximation are characterized. The approach does not involve analytic perturbation theory of linear operators and allows for multiplicities in the singular value spectrum. Our main result is some general non-asymptotic upper bound on the accuracy of approximation which involves explicitly the singular values of $C$, and which is shown to be sharp in various regimes of $C$. The results are accompanied by lower bounds under diverse  assumptions. Consequences on statistical estimation problems, in particular in the recent area of low-rank matrix recovery, are discussed.
\end{abstract}
%
%\begin{keyword}[class=AMS]
%\kwd[Primary ]{60K35}
%\kwd{60K35}
%\kwd[; secondary ]{60K35}
%\end{keyword}
%
%\begin{keyword}
%\kwd{sample}
%\kwd{\LaTeXe}
%\end{keyword}

\end{frontmatter}

\small
\section{Introduction}
As a consequence of the \cite{baiyin93} law, the maximal singular value $\lambda_{\max}(\mathrm{E})$ of an iid standard Gaussian $M\times M$-matrix $\mathrm{E}$ is equal to $2\sqrt{M}(1+o(1))$ a.s. Since in addition the sequence $\lambda_{\max}(\mathrm{E})/\sqrt{M}$ is uniformly integrable (\cite{johlin01}, Chapter 8, Theorem 2.4), the corresponding bound holds in expectation as well. Similarly, $\E\lambda_{\max}(\mathrm{E})^2=4M(1+o(1))$.  Let $\Arrowvert\cdot\Arrowvert_{S_2}$ denote the Hilbert-Schmidt or Frobenius norm. Define $\hat{\pi}_1$ to be the orthogonal projection matrix onto the one-dimensional subspace of $\R^M$ maximizing $\Arrowvert \tilde{\pi}_1\mathrm{E}\Arrowvert_{S_2}^2$ over all one-dimensional orthogonal projections $\tilde{\pi}_1$. Rewriting $\lambda_{\max}(\mathrm{E})^2=\Arrowvert\hat{\pi}_1\mathrm{E}\Arrowvert_{S_2}^2$  yields 
\begin{equation}\label{eq: intro}
\E\Arrowvert\hat{\pi}_1\mathrm{E}\Arrowvert_{S_2}^2\ =\ 4M(1+o(1)).
\end{equation}
In contrast, $\E\Arrowvert \tilde{\pi}_1\mathrm{E}\Arrowvert_{S_2}^2=M$ for every fixed $\tilde{\pi}_1$. Thus, replacing one single projection by the supremum over all projections increases the Hilbert-Schmidt norm by a positive factor:
\begin{equation}\label{eq: intro 2}
\E\Arrowvert\hat{\pi}_1\mathrm{E}\Arrowvert_{S_2}^2\ -\ \E\Arrowvert\tilde{\pi}_1\mathrm{E}\Arrowvert_{S_2}^2\ =\ 3(1+o(1))M.
\end{equation} 
This effect raises the question about the accuracy for empirical reduced-rank projections in general. 
Consider the model
\begin{equation}\label{eq: model 1}
X\ =\ C\ +\ \mathrm{E}
\end{equation} 
with a deterministic matrix $C\in\R^{M\times M}$ and $\mathrm{E}$ some centered Gaussian $M\times M$-matrix whose entries are independent with variance $\sigma^2$. Here and subsequently, let
\begin{equation}\label{eq: argmax}
\hat{\pi}_r\ :=\ \underset{~~ \tilde{\pi}_r\in\SS_{M,r}}{\mathrm{Argmax}}\Arrowvert \tilde{\pi}_rX\Arrowvert_{S_2}^2\ \ \ \text{and}\ \  \  \pi_r\ \in\ \underset{~~ \tilde{\pi}_r\in\SS_{M,r}}{\mathrm{Argmax}}\E\Arrowvert \tilde{\pi}_rX\Arrowvert_{S_2}^2%\ =\ \underset{~~ \tilde{\pi}_r\in\SS_{M,r}}{\mathrm{Argmax}}\Arrowvert \tilde{\pi}_rC\Arrowvert_{S_2}^2
\end{equation}

\vspace{-4mm}
with $\SS_{M,r}$ denoting the set of all $M\times M$-matrices representing orthogonal projections onto $r$-dimensional linear subspaces of $\R^M$. %, that is, $\SS_{M,r}$ is the Grassmann manifold, considered as a subset of $\mathbb{R}^{M\times M}$. 
How close is $\E\Arrowvert \hat{\pi}_r X\Arrowvert_{S_2}^2$ to its deterministic counterpart $\E\Arrowvert \pi_rX\Arrowvert_{S_2}^2=\Arrowvert\pi_rC\Arrowvert_{S_2}^2+\sigma^2rM$ 
if the Gaussian matrix $X$  is not centered?
For every fixed $M\in\N$ and $\sigma^2>0$, the following questions are natural:
\begin{itemize}
\item[(A)] Does there exist some favorable matrix $C=\E X$ for which the accuracy of approximation $\E\Arrowvert \hat{\pi}_rX\Arrowvert_{S_2}^2-\E\Arrowvert\pi_rX\Arrowvert_{S_2}^2$ improves over the situation described in (\ref{eq: intro 2})?
\item[(B)] Does there exist for any arbitrarily large real number $c$ some unfavorable matrix $C(c)$ such that $\E\Arrowvert \hat{\pi}_r X\Arrowvert_{S_2}^2-\E\Arrowvert \pi_rX\Arrowvert_{S_2}^2\geq c$?
\end{itemize} 
Based on the random variable $X=C+\mathrm{E}$ within model (\ref{eq: model 1}), denote the difference by  
\begin{align}
\delta_{C,M,\sigma^2,r}\ :=&\ \E\Arrowvert \hat{\pi}_r X\Arrowvert_{S_2}^2  - \E\Arrowvert \pi_rX\Arrowvert_{S_2}^2,\label{eq: delta}
\end{align}
which, in terms of singular values, is equal to
$$\sum_{i=1}^r\E\Big(\hat{\lambda}_i^2-\lambda_i^2-\sigma^2M\Big),
$$ 

\vspace{-4mm}
where $\hat{\lambda}_1\geq\hat{\lambda}_2\geq ...\geq \hat{\lambda}_M$ and  ${\lambda_1}\geq{\lambda_2}\geq ...\geq {\lambda}_M$ denote the singular values of $X$ and $C$, respectively.
The goal in the present article is to study this quantity $\delta_{C,M,\sigma^2,r}$, to derive universal upper and lower bounds,  and to characterize favorable and unfavorable types of matrices $C$ in terms of the accuracy of approximation $\delta_{C,M,\sigma^2,r}$.

The motivation for considering this problem is two-fold. First of all, as 
$$
\E\Arrowvert \hat{\pi}_rX\Arrowvert_{S_2}^2\ =\ \Big(\E\Arrowvert \hat{\pi}_r X\Arrowvert_{S_2}^2 \ -\ \E\Arrowvert \pi_rX\Arrowvert_{S_2}^2\Big)\ +\ \Arrowvert \pi_r C\Arrowvert_{S_2}^2\ +\ \sigma^2rM
$$
and $\E\Arrowvert \hat{\pi}_r X\Arrowvert_{S_2}^2  - \E\Arrowvert \pi_rX\Arrowvert_{S_2}^2=\delta_{C,M,\sigma^2,r}\geq 0$, the problem is of theoretical interest  as our results complement the bound in (\ref{eq: intro}) for centered Gaussian matrices with a detailed {\it non-asymptotic analysis} of the noncentered case, extending also to more general rank-$r$ projections. Finite-rank perturbations of random matrices  have found recently a lot of attention, see \cite{cadofe09}, \cite{cadofe12}, \cite{pireso12}, \cite{tao12} among others. \cite{tao12}, Theorem 1.7,  studies the the eigenvalue value spectrum of low rank perturbations of an iid complex random matrix and proves, as a special case, that  $\gamma_{\max}(C+\mathrm{E}/(\sigma \sqrt{M}))=\gamma_{\max}(C) + o_p(1)$ as $M\rightarrow\infty$ and $\mathrm{rank}(C)=O(1)$ as long as $\arrowvert \gamma_{\max}(C)\arrowvert =O(1)$ is sufficiently large, with $\gamma_{\max}(C)$ the eigenvalue of $C$ which is maximal in absolute value. \cite{cadofe09} and \cite{pireso12} study Wigner matrices instead of iid random matrices. Somewhat remarkably, the outlier eigenvalues of the perturbed matrix are not close in probability to those of the original matrix $C$ but to some shifted value $\lambda_i(C)+\sigma^2/\lambda_i(C)$, where $\sigma^2$ is the common variance of the entries of the Wigner matrix, and $\lambda_i(C)$ the eigenvalues of an Hermitian matrix $C$.  Our results are complementary: 
\begin{itemize}
\item We derive {\it non-asymptotic} cumulated second moment bounds on the singular values in the deformed (non-Hermitian) iid real Gaussian matrix case, i.e.~the noise level $\sigma^2$ and the dimension $M$ are fixed but arbitrary throughout the analysis, and the constants involved in our bounds do not depend on them.
\item  The perturbation matrix $C$  is not required to be of low or uniformly bounded  rank, for example, our results cover the case  $\mathrm{rank}(C)=\lfloor M/2\rfloor$ or $\mathrm{rank}(C)=M$. 
\item Our proofs differ significantly from the techniques of the above mentioned results but rely on empirical process techniques without making use of classical random matrix tools. The novelty in the proof of the subsequent Theorem \ref{thm: general} is that a slicing argument is used for bounding the expectation of the supremum over some {\it non-centered} process.
\end{itemize} 

Although our results extend without difficulties to the self-adjoint dilation $\tilde{X}$ of $X$ in $\R^{2M\times 2M}$,
 it remains open whether,  in an appropriate asymptotic sense, the eigenvalue spectrum of $\tilde{X}$ behaves similarly to  the deformed Wigner case as studied by \cite{cadofe09} and \cite{pireso12} with finite-rank perturbations, as their assumptions do not apply to this setting.

As concerns applicability in mathematical statistics, the study of $\E\Arrowvert\hat{\pi}_rX\Arrowvert_{S_2}^2=\sum_{i=1}^r\E\,\hat{\lambda}_i^2$ arises naturally when infering about quantities like
\begin{equation}\label{eq: exva}
\Arrowvert C\Arrowvert_{S_2}^2-\Arrowvert\pi_rC\Arrowvert_{S_2}^2\ \ \ \text{or}\ \ \  \arg\min_{r\geq 1}\bigg\{\frac{\Arrowvert\pi_rC\Arrowvert_{S_2}^2}{\Arrowvert C\Arrowvert_{S_2}^2}\geq \alpha\bigg\},\ \alpha\in(0,1],
\end{equation} 
which are of interest for analyzing and understanding the singular value spectrum, in particular in high dimension. As above and subsequently, for any matrix $C\in\R^{M\times M}$,  its singular values $\lambda_1,...,\lambda_M$ are ordered in decreasing magnitude. In terms of singular values,
$$
\Arrowvert C\Arrowvert_{S_2}^2\ =\ \sum_{i=1}^M\lambda_i^2\ \ \ \text{and}\ \ \ \Arrowvert \pi_rC\Arrowvert_{S_2}^2\ =\ \sum_{i=1}^r\lambda_i^2.
$$
If
$
C = \sum_{i=1}^M\lambda_iU_iV_i'
$
 denotes some singular value decomposition of $C$, where $U_1,...,U_M$ and $V_1,...,V_M$ are two sets of orthonormal vectors in $\R^M$, then the maximizer 
\begin{align*}
&\underset{~~\tilde{\pi}_r\in\SS_{M,r}}{\mathrm{Argmax}}\ \E\Arrowvert\tilde{\pi}_r(C+\mathrm{E})\Arrowvert_{S_2}^2\ =\  \underset{~~ \tilde{\pi}_r\in\SS_{M,r}}{\mathrm{Argmax}}\ \Arrowvert \tilde{\pi}_rC\Arrowvert_{S_2}^2
\end{align*} 
is unique if and only if $\lambda_r>\lambda_{r+1}$, in which case it is equal to  the orthogonal projection $\sum_{i=1}^rU_iU_i'$  onto the linear space spanned by the orthonormal  column vectors $U_1,...,U_r$, and $
\pi_rC = \sum_{i=1}^r\lambda_iU_iV_i'
$. %Note that due to potential multiplicities of singular values, the orthonormal vectors $U_1,...,U_r$ themselves in the singular value decomposition are not necessarily unique. 
%If $\lambda_r=\lambda_{r+1}=...=\lambda_{r+t}$ for some $t\leq M-r$, where $t$ is the maximal number with this property, then any orthogonal projection onto some $r$-dimensional subspace of $\mathrm{span}(U_1,...,U_{r+t})$ maximizes $\Arrowvert \tilde{\pi}_rC\Arrowvert_{S_2}$. 
In the context of covariance matrices, the ratio $\Arrowvert\pi_rC\Arrowvert_{S_2}^2/\Arrowvert C\Arrowvert_{S_2}^2$ is often referred to as percentage of the "explained variance" by the first $r$ principal components, and the second expression in (\ref{eq: exva}) determines the smallest number of principal components needed to explain a prescribed percentage $\alpha$ of the overall variance. Within our model, the statistics
$\Arrowvert \pi_rX\Arrowvert_{S_2}^2-\sigma^2rM$ estimates the expression  $\Arrowvert\pi_rC\Arrowvert_{S_2}^2$  in (\ref{eq: exva}) unbiasedly, but note that $\pi_r=\pi_r(C)$ is not available in advance as $C$ itself and in particular its singular spaces are unknown. Thus, the  first question in the analysis is whether the empirical counterpart 

\vspace{-4mm}
\begin{equation*}
\Arrowvert \hat{\pi}_rX\Arrowvert_{S_2}^2\ -\ \sigma^2rM\ =\ \sum_{i=1}^r\hat{\lambda}_i^2\ -\ \sigma^2rM
\end{equation*}

\vspace{-2mm}
does the job as well, where $\hat{\lambda}_1,...,\hat{\lambda}_r$ denote the first $r$ largest singular values of $X$. Our profound analysis about this problem will show that its bias depends strongly on the unknown matrix $C$ itself, and even for favorable rank-$r$ matrices $C=\E X$ of arbitrarily large amplitude and rank-$r$ projections, the bias  $\E\Arrowvert\hat{\pi}_r X\Arrowvert_{S_2}^2-\E\Arrowvert\pi_r X\Arrowvert_{S_2}^2$ remains of the order $\sigma^2r(M-r)$, which is shown to be unimprovable in general. Note at this point that $\E\Arrowvert\hat{\pi}_rX\Arrowvert_{S_2}^2-\Arrowvert \pi_rC\Arrowvert_{S_2}^2\geq  \E\Arrowvert{\pi}_rX\Arrowvert_{S_2}^2-\Arrowvert \pi_rC\Arrowvert_{S_2}^2=\sigma^2rM$, but there is a priori no reason why the difference 
$$\E\Arrowvert \hat{\pi}_r X\Arrowvert_{S_2}^2 \ -\ \E\Arrowvert \pi_rX\Arrowvert_{S_2}^2$$ cannot be even of substantially smaller order for "good" choices $C$, see question (A).  
In order to keep the technical expenditure as small as possible, we consider the model $X=C+\mathrm{E}$ as mentioned above, but we conjecture that similar non-asymptotic implications (somewhat different and still to be derived  for the Wishart case) will be valid for the squared Hilbert-Schmidt norm of rank-$r$ projections of high-dimensional empirical covariance matrices $YY'$, with $Y\sim\mathcal{N}(0,C)$ for some positive semidefinite matrix $C$, with consequences on the robustness properties of principal component analysis. %Some further implication of our analysis  to the concrete example of linear entropy or purity estimation in quantum information theory based on observation schemes like the ours is briefly discussed in Section 
%\ref{sec: discussion}.

 In the simplest special cases where the question is non-trivial, the main findings  of the article can be summarized as follows:
\begin{theorem}[Prototype of weak accuracy]\label{thm: intro}
Let $C_{\alpha}=\alpha\id$ for some arbitrary real number $\alpha\in\R$, where $\id$ denotes the $M\times M$-identity matrix. Then the following statements hold true:
(i) The case $C=C_{\alpha}=\alpha\id$, $\alpha\not=0$, is always worse than the case $C=0$:
$$
\delta_{C_{\alpha},M,\sigma^2,r}\ \geq\ \delta_{0,M,\sigma^2,r}.
$$
(ii) The difference $\delta_{C_{\alpha},M,\sigma^2,r}$ explodes at least linearly in the amplitude $\arrowvert\alpha\arrowvert$: There exists some constant $c_1>0$, independent of $\sigma^2$, $r$ and $M$, such that
$$ \liminf_{\arrowvert\alpha\arrowvert\rightarrow\infty}\,\frac{\delta_{C_{\alpha},M,\sigma^2,r}}{\arrowvert\alpha\arrowvert}\ \geq\ c_1\,\sigma r\sqrt{M-r}\ \ \ \text{for all}\ r\leq M-r.
$$
\end{theorem}
\begin{theorem}[Universal upper and lower bound in case $r=1$]\label{thm: universal}
(i) There exists some constant $c_2>0$ independent of $C,M,\sigma^2$, such that for all $C\in\R^{M\times M}$ and $M\geq 2$
\begin{equation}\label{eq: upper intro}
\delta_{C,M,\sigma^2,1}\ \leq\ c_2\Big(\sigma^2M+\sigma\sqrt{M}\Arrowvert C\Arrowvert_{S_{\infty}}\Big),
\end{equation}
where $\Arrowvert C\Arrowvert_{S_{\infty}}$ denotes the spectral norm of $C$.

(ii)  There exist constants $c_3>0$ and $M_0\in\N$, independent of $\sigma^2$, such that 
\begin{equation}\label{eq: intro lower}
\inf_{C\in\R^{M\times M}}\delta_{C,M,\sigma^2,1}\ \geq \ c_3\sigma^2(M-1)\ \ \ \text{for all}\ M\geq M_0.
\end{equation}
\end{theorem}
%Theorem \ref{thm: intro} (ii) provides a positive answer to question (B).
\begin{theorem}[Prototype of high accuracy]\label{thm: high} 
Let $C_{\alpha,s}=\mathrm{diag}(\alpha,,...,\alpha, 0,...,0)$ with $\mathrm{rank}(C_{\alpha})=s$. Then the following statements hold true:
(i) There exists some constant $c_4>0$ independent of $\alpha,r,M$ and $\sigma^2$, such that
\begin{equation}\label{eq: intro 1}
\delta_{C_{\alpha,r},M,\sigma^2,r}\ \leq\ c_4\,\sigma^2rM\ \ \ \text{for $r\leq M-r$ and \underline{every} $\alpha\in\R$.}
\end{equation}
(ii) There exists some constant $c_5>0$ independent of $M$, $r$ and $\sigma^2$, such that  the bound (\ref{eq: intro}) is asymptotically sharp  in the following sense:
$$
\liminf_{\arrowvert\alpha\arrowvert\rightarrow\infty}\max_{s\in\{r,M-r\}}\delta_{C_{\alpha,s},M,\sigma^2,s}\ \geq\ c_5\sigma^2r(M-r).
$$
The same result holds true even without the maximum over $\{1,M-1\}$ for $s=r=1$.
\end{theorem}

\vspace{-1.5mm}
Theorem \ref{thm: high} (i) describes some special case of the more general upper bound in Theorem \ref{thm: general}, which applies for every matrix $C\in\R^{M\times M}$ of $\mathrm{rank}(C)\geq r$. In case that $\lambda_i=\alpha$ for $i\leq r$ and $\lambda_i=\beta<\alpha$ for $i> r$, the bound of Theorem \ref{thm: general} approaches the unimprovable upper bound for matrices $C$ of constant singular value spectrum as $\beta\rightarrow\alpha$.

%The study of properties of the singular value spectrum of iid matrix ensembles has a long history in various fields of mathematics which we cannot satisfactorily summarize at this point. Besides the convergence properties of the empirical spectral measure to the limiting spectral distribution, the behavior of the largest singular value of centered matrices with independent entries and growing dimension has been studied in detail (see, for instance, the monograph of Guionnet (2008) for comprehensive literature on this topic). For the specific case $r=1$, our second motivation is to complement the bound in expectation (\ref{eq: intro}) with a detailed study for the non-centered case, discovering the strong dependence on the mean $C=\E X$ in such a bound. Moreover, the results are generalizing also to arbitrary rank-$r$ projections for $r>1$.

The article is organized as follows. In Section \ref{sec: pre}, we introduce the notation and describe some basic observations about $\delta_{C,M,\sigma^2,r}$. Prototypes of matrices $C$ of weak accuracy and first lower bounds are studied in Section \ref{sec: weak accuracy}. The supremum over the centered differences
$$
\sup_{\tilde{\pi}_r\in\SS_{M,r}}\Big(\Arrowvert \tilde{\pi}_r(C+\mathrm{E})\Arrowvert_{S_2}^2-\E\Arrowvert \tilde{\pi}_r(C+\mathrm{E})\Arrowvert_{S_2}^2\ -\ \Big[\Arrowvert \pi_r(C+\mathrm{E})\Arrowvert_{S_2}^2-\E\Arrowvert \pi_r(C+\mathrm{E})\Arrowvert_{S_2}^2\Big]\Big)
$$
is analyzed in Section \ref{sec: initial}. The process of centered differences and modifications thereof are central for our analysis. Our main results are given in Section \ref{subsec: iteration}. Our general idea on how to derive potentially sharp upper bounds on $\delta_{C,M,\sigma^2,r}$ for general $M\times M$-matrices $C$ is described at the beginning of that section.  The upper bounds are complemented with  lower bounds in Section \ref{sec: lower bounds}. %More general $L^p$-bounds are presented for rank-$r$-matrices as well. As a statistical application of our results, we propose an estimator for $\Arrowvert C\Arrowvert_{S_2}^2$ and prove its asymptotic normality and efficiency in a high-dimensional setting when $C$ is of low rank approximately. 
Consequences on statistical estimation problems, in particular in the recent area of low-rank matrix recovery,  are discussed in Section \ref{sec: discussion}. Section \ref{sec: theorem main} is devoted to the proof of Theorem \ref{thm: general}. The proof of Theorem \ref{thm: lower bounds} is deferred to Section \ref{sec: lower bounds 3}.

\vspace{-2mm}
\section{Preliminaries}\label{sec: pre}

\vspace{-1mm}
\subsection{Notation}The notation $\lesssim$ means less or equal up to some non-negative multiplicative constant which does not depend on the variable parameters in the expression. $A\sim B$ should be read as $A\lesssim B$ and $B\lesssim A$ at once.
If not stated otherwise, $\mathrm{E}$ is a centered Gaussian matrix whose entries are independent with variance $\sigma^2$.  Subsequently, $\Arrowvert .\Arrowvert_{S_p}$, $1\leq p\leq\infty$, denotes the Schatten-$p$-norm on $\R^{M\times M}$, i.e.~for any $C\in\R^{M\times M}$, the $\Arrowvert C\Arrowvert_{S_p}$ coincides with the $\ell_p$-norm of its singular values $\lambda_1\geq \lambda_2\geq ...\geq\lambda_M$:

\vspace{-4mm}
$$
\Arrowvert C\Arrowvert_{S_p}\ =\ \bigg(\sum_{i=1}^M\lambda_i^p\bigg)^{1/p}\ \ \text{for $1\leq p<\infty$,\ \  and}\ \ \ \Arrowvert C\Arrowvert_{S_{\infty}}\ =\ \lambda_1.
$$
Specifically, $\Arrowvert\cdot\Arrowvert_{S_1}$, $\Arrowvert\cdot\Arrowvert_{S_2}$ and $\Arrowvert \cdot\Arrowvert_{S_{\infty}}$ are referred to as nuclear norm, Hilbert-Schmidt or Frobenius norm, and spectral norm, respectively. $\tr(C)$ denotes the trace of $C\in\R^{M\times M}$. For any $A\in\R^{d\times M}$, $A'$ denotes its transpose in $\R^{M\times d}$. $\id$ denotes the $M\times M$-identity matrix, and $\id_r$ the diagonal matrix $\mathrm{diag}(1,...,1,0,...,0)$ of rank $r$. As usual, $O(M)$ describes the orthogonal group, i.e.~the group  of orthogonal $M\times M$-matrices. For any totally-bounded, pseudometric space $(\mathcal{X},d)$ and any subset $E\subset\XX$, the covering number $N(E,d,\delta)$ is the smallest number of closed $d$-balls in $\XX$ of radius $\delta$ needed to cover $E$.

\vspace{-2mm}
\subsection{Some basic observation about $\delta_{C,M,\sigma^2,r}$}
The following representation clarifies the problem under consideration.
For some arbitrary matrix $A\in\R^{M\times M}$, $\Arrowvert A\Arrowvert_{S_{\infty}}=1$, and $\alpha\in\R$, inspection of the quantity $\delta_{\alpha A,\sigma^2,M,r}$ shows
\begin{align*}
\delta_{\alpha A,M,\sigma^2,r}\ &=\ \E\big\Arrowvert \hat{\pi}_r (\alpha A+\mathrm{E})\big\Arrowvert_{S_2}^2 \ -\ \E\big\Arrowvert \pi_r ( \alpha A+\mathrm{E})\big\Arrowvert_{S_2}^2\\ 
&=\ \E\bigg(\sup_{\tilde{\pi}_r\in\SS_{M,r}}\bigg\{\big\Arrowvert\tilde{\pi}_r\mathrm{E}\big\Arrowvert_{S_2}^2-\big\Arrowvert{\pi}_r\mathrm{E}\big\Arrowvert_{S_2}^2 +2\alpha\, \tr\Big(\mathrm{E}'(\tilde{\pi}_r-\pi_r)A\Big)\\
&\ \ \ \ \ \ \ \ \ \ \ \ \ \ \ \ \ \ \ \ \ \ \ \ \ \ \ \ \ \ \ \ \ \ \ \ \ \ \ \ \ \ \ \ \ \ \ \ \underbrace{-\alpha^2\Big(\big\Arrowvert{\pi}_rA\big\Arrowvert_{S_2}^2-\big\Arrowvert \tilde{\pi}_rA\big\Arrowvert_{S_2}^2\Big)}_{\textrm{"compensation term"},\ \leq 0}\bigg\}\bigg).
\end{align*}
The two processes $$
\Big(\Arrowvert\tilde{\pi}_r\mathrm{E}\Arrowvert_{S_2}^2-\Arrowvert{\pi}_r\mathrm{E}\Arrowvert_{S_2}^2\Big)_{\tilde{\pi}_r\in\SS_{M,r}}\ \ \ \text{and}\ \ \ \Big(\alpha\tr\big(\mathrm{E}'(\tilde{\pi}_r-\pi_r)A\big)\Big)_{\tilde{\pi}_r\in\SS_{M,r}}
$$ 
are centered, while the deterministic compensation term $-\alpha^2\big(\Arrowvert{\pi}_rA\Arrowvert_{S_2}^2-\Arrowvert \tilde{\pi}_rA\Arrowvert_{S_2}^2\big)$ is less or equal to zero, for any choice of $A$ and every $\tilde{\pi}_r\in\SS_{M,r}$. Note that the stochastic term  $\alpha\tr\big(\mathrm{E}'(\tilde{\pi}_r-\pi_r)A\big)$ is linear in $\alpha$, but the deterministic compensation term depends quadratically on  $\alpha$. This representation suggests that an interplay of amplitude $\arrowvert \alpha\arrowvert$ and structure of $A$ determine the accuracy of approximation.

\section{The prototype of weak accuracy -- the case without deterministic compensation term}\label{sec: weak accuracy}
%In principle, an asymptotically sharp lower bound of the expression 
%$
%\E\big(\sup_{~\tilde{\pi}_r\in\SS_{M,r}}\Arrowvert \tilde{\pi}_r\mathrm{E}\Arrowvert_{S_2}^2\big) - \sigma^2rM
%$ 
%can be deduced asymptotically via the semi-circular law.  

Recall the definition (\ref{eq: delta}), with $\hat{\pi}_r$ and $\pi_r$ as defined in (\ref{eq: argmax}). The first result is a lower bound on the expected squared Hilbert-Schmidt norm of the rank-$r$-projection in case that the singular value spectrum of $C$ is constant.
\begin{proposition}\label{prop: stability}
Let $C_{\alpha}\in\R^{M\times M}$ with singular value decomposition $U\Lambda_{\alpha} V'$. Assume that $\Lambda_{\alpha}=\alpha \id$ with  some non-negative number $\alpha\in\R$.  Then 

(i) $\ \ \ \delta_{C_{\alpha},M,\sigma^2,r}\ \geq\ \delta_{0,M,\sigma^2,r}\,$ for every $\alpha>0$, 

and

(ii)
$\ \ \  \displaystyle
\liminf_{ \alpha\rightarrow\infty}\,\frac{\delta_{C_{\alpha},M,\sigma^2,r}}{\alpha}\,
\gtrsim\ r\sigma\sqrt{M-r}
$ \ for any $r\leq M-r$.
\end{proposition}

\begin{remark}
Proposition \ref{prop: stability} (i) demonstrates that the accuracy in case $C=\alpha\id$ is always worse than in case $C=0$. (ii) complements this observation with an asymptotic lower bound: the difference  $\delta_{C_{\alpha},M,\sigma^2,r}$ explodes at least linearly in the amplitude $\alpha$. In particular, (ii) provides a positive answer to question (B) in the Introduction. 
\end{remark}

%Proposition \ref{prop: stability} indicates that stability of the reduced-rank projection is coupled not on the amplitude of the singular value spectrum only, and the direction of the inequality in (i)  may be surprising at first glance. %As demonstrated in (ii), depending on the geometry of $C$, the concentration may get even weaker the larger the amplitude $\Arrowvert C\Arrowvert_{S_{\infty}}$ is. 

\vspace{-3mm}
\paragraph{\sc Proof}$Y=_{\DD}Z$ for two random variables $Y$ and $Z$ means that their distributions coincide. Since $U'\tilde{\pi}_rU\in\SS_{M,r}$ for any $\tilde{\pi}_r\in\SS_{M,r}$, $U'\mathrm{E} V=_{\DD}\mathrm{E}$ for two fixed orthogonal matrices $U$ and $V$, and $\Arrowvert U'AV\Arrowvert_{S_2}^2=\Arrowvert A\Arrowvert_{S_2}^2$ for any $A\in\R^{M\times M}$, we may assume without loss of generality that $C=\alpha \id$. Let $\pi_r$ be as given in (\ref{eq: argmax}), i.e.~in this case, $\pi_r$ denotes some arbitrary but fixed element of $\SS_{M,r}$. First, 
$$
\big\Arrowvert\tilde{\pi}_r(\alpha  \id+\mathrm{E})\big\Arrowvert_{S_2}^2-\big\Arrowvert{\pi}_r(\alpha \id+\mathrm{E})\big\Arrowvert_{S_2}^2\ =\ 2\alpha\tr\Big(\mathrm{E}'\big(\tilde{\pi}_r-\pi_r\big)\Big)\ +\ \Arrowvert \tilde{\pi}_r\mathrm{E}\Arrowvert_{S_2}^2-\Arrowvert \pi_r\mathrm{E}\Arrowvert_{S_2}^2,
$$
i.e. we need a lower bound on
$$
\E\left(\sup_{\tilde{\pi}_r\in\SS_{M,r}}2\alpha\tr\Big(\mathrm{E}'\big(\tilde{\pi}_r-\pi_r\big)\Big)\ +\ \Arrowvert \tilde{\pi}_r\mathrm{E}\Arrowvert_{S_2}^2-\Arrowvert \pi_r\mathrm{E}\Arrowvert_{S_2}^2\right).
$$
Note that the supremum within the expectation is non-negative simply because $\pi_r\in\SS_{M,r}$. Let $\mathrm{E}^*\in\R^{M\times M}$ be the matrix with the entries
$
\mathrm{E}_{ij}^* = -\mathrm{E}_{ij}$, $i,j=1,...,M
$. 
 By the symmetry of the Gaussian distribution,
\begin{align*}
2\alpha\tr\Big(\mathrm{E}'\big(\tilde{\pi}_r-&\pi_r\big)\Big)+ \Arrowvert \tilde{\pi}_r\mathrm{E}\Arrowvert_{S_2}^2-\Arrowvert \pi_r\mathrm{E}\Arrowvert_{S_2}^2\\
& =_{\DD}\ 2\alpha\tr\Big({\mathrm{E}^*}'\big(\tilde{\pi}_r-\pi_r\big)\Big) + \Arrowvert \tilde{\pi}_r\mathrm{E}^*\Arrowvert_{S_2}^2-\Arrowvert \pi_r\mathrm{E}^*\Arrowvert_{S_2}^2\\
&=\ -2\alpha\tr\Big(\mathrm{E}'\big(\tilde{\pi}_r-\pi_r\big)\Big) + \Arrowvert \tilde{\pi}_r\mathrm{E}\Arrowvert_{S_2}^2-\Arrowvert \pi_r\mathrm{E}\Arrowvert_{S_2}^2.
\end{align*}
Consequently,
\begin{align}
\E\bigg(&\sup_{\tilde{\pi}_r\in\SS_{M,r}}2\alpha\tr\Big(\mathrm{E}'\big(\tilde{\pi}_r-\pi_r\big)\Big)\ +\ \Arrowvert \tilde{\pi}_r\mathrm{E}\Arrowvert_{S_2}^2-\Arrowvert \pi_r\mathrm{E}\Arrowvert_{S_2}^2\bigg)\nonumber\\
&=\ \frac{1}{2} \E\biggl(\sup_{\tilde{\pi}_r\in\SS_{M,r}}2\alpha\tr\Big(\mathrm{E}'\big(\tilde{\pi}_r-\pi_r\big)\Big)\ +\ \Arrowvert \tilde{\pi}_r\mathrm{E}\Arrowvert_{S_2}^2-\Arrowvert \pi_r\mathrm{E}\Arrowvert_{S_2}^2\biggr)\nonumber\\
&\ \ \ \ \ \ +\  \frac{1}{2}\E\biggl(\sup_{\tilde{\pi}_r\in\SS_{M,r}}2\alpha\tr\Big({\mathrm{E}^*}'\big(\tilde{\pi}_r-\pi_r\big)\Big)\ +\ \Arrowvert \tilde{\pi}_r\mathrm{E}^*\Arrowvert_{S_2}^2-\Arrowvert \pi_r\mathrm{E}^*\Arrowvert_{S_2}^2\biggr)\nonumber\\
&=\ \frac{1}{2}\E\bigg[\sup_{\tilde{\pi}_r\in\SS_{M,r}}2\alpha\tr\Big(\mathrm{E}'\big(\tilde{\pi}_r-\pi_r\big)\Big)\ +\ \Arrowvert \tilde{\pi}_r\mathrm{E}\Arrowvert_{S_2}^2-\Arrowvert \pi_r\mathrm{E}\Arrowvert_{S_2}^2\nonumber\\
&\ \ \ \ \ \ \ \ \ \ \ \ +\  \sup_{\tilde{\pi}_r\in\SS_{M,r}}2\alpha\tr\Big({\mathrm{E}^*}'\big(\tilde{\pi}_r-\pi_r\big)\Big)\ +\ \Arrowvert \tilde{\pi}_r\mathrm{E}^*\Arrowvert_{S_2}^2-\Arrowvert \pi_r\mathrm{E}^*\Arrowvert_{S_2}^2\bigg]\nonumber\\
&\geq\  \frac{1}{2}\E\bigg[\sup_{\tilde{\pi}_r\in\SS_{M,r}}\bigg(2\alpha\tr\Big(\mathrm{E}'\big(\tilde{\pi}_r-\pi_r\big)\Big)\ +\ \Arrowvert \tilde{\pi}_r\mathrm{E}\Arrowvert_{S_2}^2-\Arrowvert \pi_r\mathrm{E}\Arrowvert_{S_2}^2\nonumber\\
&\ \ \ \ \ \ \ \ \ \ \ \ \ \ \ \ \ \ \ \ \ \ \ \ +\  2\alpha\tr\Big({\mathrm{E}^*}'\big(\tilde{\pi}_r-\pi_r\big)\Big)\ +\ \Arrowvert \tilde{\pi}_r\mathrm{E}^*\Arrowvert_{S_2}^2-\Arrowvert \pi_r\mathrm{E}^*\Arrowvert_{S_2}^2\bigg)\bigg]\nonumber\\
&=\ \E\bigg(\sup_{\tilde{\pi}_r\in\SS_{M,r}}\Arrowvert \tilde{\pi}_r\mathrm{E}\Arrowvert_{S_2}^2-\Arrowvert \pi_r\mathrm{E}\Arrowvert_{S_2}^2\bigg),\ \ \ \ \ \ \ \ \ \ \  \nonumber
%\ =\ \E\bigg(\sup_{\tilde{\pi}_r\in\SS_{M,r}}\Arrowvert \tilde{\pi}_r\mathrm{E}\Arrowvert_{S_2}^2\bigg)\ -\ \sigma^2rM,\nonumber
\end{align}
which proves part (i) of the proposition. As concerns the proof of (ii), observe that %for any $\alpha>0$
\begin{align}
\E\biggl(\sup_{\tilde{\pi}_r\in\SS_{M,r}}&2\alpha\tr\Big(\mathrm{E}'\big(\tilde{\pi}_r-\pi_r\big)\Big)\ +\ \Arrowvert \tilde{\pi}_r\mathrm{E}\Arrowvert_{S_2}^2-\Arrowvert \pi_r\mathrm{E}\Arrowvert_{S_2}^2\biggr)\nonumber\\
&\geq\ 2\alpha\,\E\biggl(\sup_{\tilde{\pi}_r\in\SS_{M,r}}\tr\Big(\mathrm{E}'\big(\tilde{\pi}_r-\pi_r\big)\Big)\biggr)\ - \ \E\biggl(\sup_{\tilde{\pi}_r\in\SS_{M,r}}\Arrowvert {\pi}_r\mathrm{E}\Arrowvert_{S_2}^2-\Arrowvert \tilde{\pi}_r\mathrm{E}\Arrowvert_{S_2}^2\biggr)\nonumber\\
&\geq\ 2\alpha\,\E\biggl(\sup_{\tilde{\pi}_r\in\SS_{M,r}}\tr\Big(\mathrm{E}'\big(\tilde{\pi}_r-\pi_r\big)\Big)\biggr)\ -\ \sigma^2rM,\label{eq: alpha}
\end{align}
since 
$\E\big(\sup_{\tilde{\pi}_r\in\SS_{M,r}}\Arrowvert {\pi}_r\mathrm{E}\Arrowvert_{S_2}^2-\Arrowvert \tilde{\pi}_r\mathrm{E}\Arrowvert_{S_2}^2\big)\ \leq\ \sigma^2rM$. By Sudakov's minoration, there exists some universal constant $c_{\textrm{Sud}}$ such that
\begin{equation}\label{eq: sud1}
\E\sup_{\tilde{\pi}_r\in\SS_{M,r}}\tr\Big(\mathrm{E}'\big(\tilde{\pi}_r-\pi_r\big)\Big)\ \geq \ \sigma\,c_{\textrm{Sud}}\,\delta\sqrt{\log N\big(\SS_{m,r}, d_{S_2},\delta\big)}
\end{equation}
for any $\delta>0$. Proposition 8 in \cite{pajo98} states for any $r\leq M-r$
$$
\Big(\frac{c'}{\xi}\Big)^{r(M-r)}\ \leq\ N\big(\SS_{M,r}, d_{S_2},\xi \sqrt{r}\big),\ \ \ \forall\ \xi>0,
$$
with some universal constant $c'>0$ which does not depend on $r$ and $M$. Choosing $\delta=\sqrt{r}c'/\e$ and plugging (\ref{eq: sud1}) into (\ref{eq: alpha}) yields, for some constant $c$ which does not depend on $\sigma^2,r,M$ and $\alpha$: 
$$
\delta_{\alpha\id,M,\sigma^2,r}\ \geq\ c\sigma r\alpha\sqrt{M-r}-\sigma^2rM.
$$
Dividing both sides by $\alpha$ and taking the limes inferior proves (ii). \hfill$\square$

In view of the representation in Section \ref{sec: pre}, the case $C=\alpha\id$ is the prototype of weak accuracy as there is no deterministic compensation term in the expression of the supremum. It follows from the subsequent Corollary \ref{cor: upper} that the lower bound of Proposition \ref{prop: stability} (ii) is sharp. 

%Using that $(1/2)a^p+(1/2)b^p\geq c_p(a+b)^p$ for any $1\leq p<\infty$, $a,b\geq 0$, and some real non-negative constant $c_p$, we obtain similarly 
%$$
%\left\Arrowvert\sup_{\tilde{\pi}_r\in\SS_{M,r}}\big\Arrowvert\tilde{\pi}_r(C_{\alpha}+\mathrm{E})\big\Arrowvert_{S_2}^2-\big\Arrowvert{\pi}_r(C_{\alpha}+\mathrm{E})\big\Arrowvert_{S_2}^2\right\Arrowvert_{L_p} \gtrsim \bigg\Arrowvert\sup_{\tilde{\pi}_r\in\SS_{M,r}}\Arrowvert \tilde{\pi}_r\mathrm{E}\Arrowvert_{S_2}^2-\Arrowvert \pi_r\mathrm{E}\Arrowvert_{S_2}^2\bigg\Arrowvert_{L_p}.
%$$
%As a hint towards the influence of structure and amplitude for "good matrices" 

\section{$S_2$-$S_{\infty}$-chaining  bounds for the supremum over the centered process and first consequences on $\delta_{C,M,\sigma^2,r}$}\label{sec: initial}
Let $X=C+\mathrm{E}$ as described in (\ref{eq: model 1}). Recall  the definition from the Introduction
$$
\pi_r\in\underset{~~\tilde{\pi}_r\in\SS_{M,r}}{\mathrm{Argmax}}\E\Arrowvert\, \tilde{\pi}_rX\Arrowvert_{S_2}^2\ = \ \underset{~~\tilde{\pi}_r\in\SS_{M,r}}{\mathrm{Argmax}}\Arrowvert\, \tilde{\pi}_rC\Arrowvert_{S_2}^2.
$$ 
Because of $\E\Arrowvert \pi_r(C+\mathrm{E})\Arrowvert_{S_2}^2\geq \E\Arrowvert\tilde{\pi}_r(C+\mathrm{E})\Arrowvert_{S_2}^2$ for every $\tilde{\pi}_r\in\SS_{M,r}$, it follows that
\begin{align}
\E\sup_{\tilde{\pi}_r\in\SS_{M,r}}&\Big(\Arrowvert \tilde{\pi}_r(C+\mathrm{E})\Arrowvert_{S_2}^2-\Arrowvert \pi_r(C+\mathrm{E})\Arrowvert_{S_2}^2\Big)\nonumber\\
&\leq\
\E\sup_{\tilde{\pi}_r\in\SS_{M,r}}\Big(\Arrowvert \tilde{\pi}_r(C+\mathrm{E})\Arrowvert_{S_2}^2-\E\Arrowvert \tilde{\pi}_r(C+\mathrm{E})\Arrowvert_{S_2}^2\label{eq: zent1}\\ 
&\ \ \ \ \ \ \ \ \  \ \ \ \ \ \ \ \ \ \ -\ \Big[\Arrowvert \pi_r(C+\mathrm{E})\Arrowvert_{S_2}^2-\E\Arrowvert \pi_r(C+\mathrm{E})\Arrowvert_{S_2}^2\Big]\Big).\nonumber
\end{align}
That is, the study of the supremum over the centered process 
$$
Z\ :=\ \sup_{\tilde{\pi}_r\in\SS_{M,r}}\Big(\Arrowvert \tilde{\pi}_r(C+\mathrm{E})\Arrowvert_{S_2}^2-\E\Arrowvert \tilde{\pi}_r(C+\mathrm{E})\Arrowvert_{S_2}^2\ -\ \Big[\Arrowvert \pi_r(C+\mathrm{E})\Arrowvert_{S_2}^2-\E\Arrowvert \pi_r(C+\mathrm{E})\Arrowvert_{S_2}^2\Big]\Big).
$$
yields some first (possibly very rough) estimate of $\delta_{C,M,\sigma^2,r}$ from above. Variants thereof are central  for our subsequent analysis in Section \ref{subsec: iteration}.%, we provide in this Section some refined study by deriving an exponential tail inequality for $Z$.
%
%There exists by now an extensive literature on large deviation inequalities for suprema of stochastic processes. For Gaussian processes, \cite{suci74}, \cite{bor75} derived a subgaussian concentration inequality of the supremum around its mean.  Gaussian concentration equally extends to $1$-Lipschitz functions, but this argument is not applicable here. Doubtless, the most famous result is due to  \cite{tala95} who established a Bernstein-type version for  bounded empirical processes. The closest to ours is Theorem 5.1 of \cite{bara10}, adapted to  random matrices and the Grassmann manifold.
%
%\vspace{-4mm}
\begin{proposition}\label{prop: rough bound}
Let $(\mathrm{E}_{ij})_{i,j=1}^M$ be a centered matrix of iid Gaussian entries with variance $\sigma^2$. Then there exists some constant $c>0$ such that for every $1\leq r<M\in\N$ and every $C\in\R^{M\times M}$
\begin{equation}\label{eq: concentration}
\E\,Z\ \leq \ c\Big(\sigma^2rM +  \sigma r\sqrt{M}\Arrowvert C\Arrowvert_{S_{\infty}}\Big).
\end{equation}
\end{proposition}
The proof of Proposition \ref{prop: rough bound} is deferred to the end of this section. We draw some first consequences on $\delta_{C,M,\sigma^2,r}$.% see (\ref{eq: upper intro})

\begin{corollary}\label{cor: upper} 
(i) (Universal upper bound) For all $C\in\R^{M\times M}$,
\begin{equation*}
\delta_{C,M,\sigma^2,r}\ \lesssim\ \sigma^2rM+\sigma r\sqrt{M}\Arrowvert C\Arrowvert_{S_{\infty}}.
\end{equation*}
(ii) (Upper bound in the "small amplitude" regime) 
\begin{align*}
\sup_{\substack{C\in\R^{M\times M}:\\\Arrowvert C\Arrowvert_{S_{\infty}}\leq \sigma\sqrt{M}}}\, \delta_{C,M,\sigma^2,r}\  \lesssim\ \sigma^2rM.
\end{align*} 
\end{corollary}
{\sc Proof} of the Corollary \ref{cor: upper} (i) follows from (\ref{eq: zent1}) and the bound on $\E\,Z$ given in Proposition \ref{prop: rough bound}; (ii) follows from (i).
\begin{remark}
It is worth being mentioned that $\E\,Z$ grows linearly with $\Arrowvert C\Arrowvert_{S_{\infty}}$, and this linear dependence is optimal for matrices $X=\alpha\id +\mathrm{E}$, $\alpha\in\R$, in view of Proposition \ref{prop: stability}. In particular, the lower bound of Proposition \ref{prop: stability} (ii) is sharp, and the universal upper bound cannot be improved without further structural assumptions on $C$.
\end{remark}
The next lemma complements the upper bound of the Corollary \ref{cor: upper} (ii) in the small amplitude regime with a lower bound.
\begin{lemma}\label{prop: stability2} 
For any real constant $\kappa>0$,
\begin{align*}
\inf_{\substack{C\in\R^{M\times M}:\\ \Arrowvert C\Arrowvert_{S_{\infty}}\,\leq\, \kappa \sigma  \sqrt{M}}}\delta_{C,M,\sigma^2,r}\ &\geq\ \delta_{0,M,\sigma^2,r}\ - \kappa^2\sigma^2rM.
\end{align*}
\end{lemma}
\paragraph{\sc Proof}With the same symmetry argument as used in the proof of Proposition \ref{prop: stability} (i), we obtain for any $C\in\R^{M\times M}$ with $\Arrowvert C\Arrowvert_{S_{\infty}}\leq\kappa \sigma\sqrt{M}$
\begin{align}
\E\biggl(\sup_{\tilde{\pi}_r\in\SS_{M,r}}&\big\Arrowvert\tilde{\pi}_r(C+\mathrm{E})\big\Arrowvert_{S_2}^2-\big\Arrowvert{\pi}_r(C+\mathrm{E})\big\Arrowvert_{S_2}^2\biggr)\nonumber\\
&\geq\ \E\bigg(\sup_{\tilde{\pi}_r\in\SS_{M,r}}\Arrowvert \tilde{\pi}_r\mathrm{E}\Arrowvert_{S_2}^2-\Arrowvert \pi_r\mathrm{E}\Arrowvert_{S_2}^2+\Arrowvert\tilde{\pi}_rC\Arrowvert_{S_2}^2-\Arrowvert\pi_rC\Arrowvert_{S_2}^2\bigg)\nonumber\\
&\geq \ \E\bigg(\sup_{\tilde{\pi}_r\in\SS_{M,r}}\Arrowvert \tilde{\pi}_r\mathrm{E}\Arrowvert_{S_2}^2-\Arrowvert \pi_r\mathrm{E}\Arrowvert_{S_2}^2\bigg) -\kappa^2\sigma^2rM.\nonumber%\\
%&\gtrsim\ \E\bigg(\sup_{\tilde{\pi}_r\in\SS_{M,r}}\Arrowvert \tilde{\pi}_r\mathrm{E}\Arrowvert_{S_2}^2-\Arrowvert \pi_r\mathrm{E}\Arrowvert_{S_2}^2\bigg)\nonumber
\end{align}

\vspace{-10mm}
$~$\hfill $\square$
%In particular for $\mathrm{rank}(C)=r=1$, the bound is of the order $\sigma^2rM$ for $\kappa=1$.

\medskip
\begin{remark}
Together with the upper Bai-Yin bound in expectation (\ref{eq: intro}), Lemma \ref{prop: stability2} implies in particular for $r=1$ that there exists some $M_0\in\N$, independent of $\sigma^2$, such that
\begin{align*}
\inf_{\substack{C\in\R^{M\times M}:\\ \Arrowvert C\Arrowvert_{S_{\infty}}\,\leq\,  \sigma  \sqrt{M}}}\delta_{C,M,\sigma^2,1}\ &\geq\ \frac{1}{2}\,\delta_{0,M,\sigma^2,1}\ \sim\ \sigma^2M\ \ \ \text{for all}\ M\geq M_0.
\end{align*}
That is, in the small amplitude regime $\Arrowvert C\Arrowvert_{S_{\infty}}\leq\sigma\sqrt{M}$, the accuracy is never better in order than in the case $C=0$, independent of the specific structure of the matrix $C$. 
\end{remark}

\paragraph{\sc Proof of Proposition \ref{prop: rough bound}} 
For any $r$-dimensional subspace $U\subset \R^M$, let $P_U\in\SS_{M,r}$ denote the orthogonal projection onto $U$. The proof is based on the classical generic chaining device. In order to make this  technique applicable, we need to investigate pairwise differences of the centered process $Z^{\sigma,M,r}$ which is pointwise given by
\begin{align*}
Z_{P_U}^{\sigma,M,r} :=&\ \Big(\Arrowvert P_U(C+\mathrm{E})\Arrowvert_{S_2}^2-\E\Arrowvert P_U(C+\mathrm{E})\Arrowvert_{S_2}^2\Big)\\
=&\ \Big(\tr\big(\mathrm{E}'P_U\mathrm{E}\big) - rM\sigma^2 +2\tr\big(C'P_U\mathrm{E}\big)\Big),\ \ P_U\in\SS_{M,r}.
\end{align*}
Denote $\pi_r^{(1)}=P_{U_1}\in\SS_{M,r}$, $\pi_r^{(2)}=P_{U_2}\in\SS_{M,r}$, and $A=\pi_r^{(1)}-\pi_r^{(2)}$. Recall that  $P_{U_i}'=P_{U_i}$ and $P_{U_i}^2=P_{U_i}$ for $i=1,2$. For any $B=(b_{ij})_{i,j=1}^M\in\R^{M\times M}$, $\mathrm{vec}(B)$ denotes the associated vector obtained by sticking together its columns,
$$\mathrm{vec}(B)\ :=\ (b_{11},\cdots,b_{M1},b_{12},\cdots,b_{M2},\cdots,b_{MM})'\in\R^{M^2}. 
$$
Observe furthermore that
\begin{align*}
\Arrowvert P_{U_1}\mathrm{E}\Arrowvert_{S_2}^2 \ -\ \Arrowvert P_{U_2}\mathrm{E}\Arrowvert_{S_2}^2\
&=\ \tr\big(\mathrm{E}'\big(P_{U_1}'P_{U_1}-P_{U_2}'P_{U_2}\big)\mathrm{E}\big)\\
%=\ \tr\big(\mathrm{E}'\big(P_{U_1}-P_{U_2}\big)\mathrm{E}\big)\\
&=\  \sum_{l,k,i=1}^MA_{lk}\mathrm{E}_{ki}\mathrm{E}_{li}\ =\ \mathrm{vec}(\mathrm{E})' \tilde{A}\,\mathrm{vec}(\mathrm{E}),
\end{align*}
where $\tilde{A}$ denotes the block-diagonal matrix $\mathrm{diag}(A,...,A)\in\R^{M^2\times M^2}$,
and, analogously,
$$
\tr\big(C'P_{U_1}\mathrm{E}\big)-\tr\big(C'P_{U_2}\mathrm{E}\big)\ =\ \mathrm{vec}(C)'\tilde{A}\,\mathrm{vec}(\mathrm{E}).
$$
Noting that $\Arrowvert \tilde{A}\Arrowvert_{S_2}=\sqrt{M}\Arrowvert A\Arrowvert_{S_2}$ and $\Arrowvert \tilde{A}\Arrowvert_{S_{\infty}}=\Arrowvert A\Arrowvert_{S_{\infty}}$, Bernstein's  inequality for quadratic forms of Gaussian variables (see, for instance, \cite{bec09}, Lemma 0.2) yields the exponential bound
\begin{align}
\mathbb{P}\Bigg( Z_{P_{U_1}}^{\sigma,M,,r}- Z_{P_{U_2}}^{\sigma,M,r} \geq\ \ \ \ \ \ \ \ \ \  \ \  \ \ \   \  \ \ \  \ \ \ \ \ \ \ \  \ \ \ \ &\nonumber\\  
2\sqrt{\sigma^4M\Arrowvert P_{U_1}-P_{U_2}\Arrowvert_{S_2}^2 +2\sigma^2\Arrowvert \tilde{A}\mathrm{vec}(C)\Arrowvert_2^2}&
\cdot\sqrt{t} + 2\sigma^2\Arrowvert P_{U_1}-P_{U_2}\Arrowvert_{S_{\infty}}t\Bigg)\ %\\ 
%&\quad\quad\quad\quad\quad\quad\quad
\leq\ \exp(-t)\label{eq: Bernstein}
\end{align}
for all $t>0$. Note that the bound is fully symmetric in $U_1$ and $U_2$.
Since 
$
\Arrowvert \tilde{A}\mathrm{vec}(C)\Arrowvert_2^2 = \Arrowvert AC\Arrowvert_{S_2}^2 \leq \Arrowvert A'A\Arrowvert_{S_1}\Arrowvert C'C\Arrowvert_{S_{\infty}}  = \Arrowvert A\Arrowvert_{S_2}^2\Arrowvert C\Arrowvert_{S_{\infty}}^2
$,
it follows that
\begin{align}
&\sqrt{\sigma^4M\Arrowvert P_{U_1}-P_{U_2}\Arrowvert_{S_2}^2+2{\sigma^2}\Arrowvert \tilde{A}\mathrm{vec}(C)\Arrowvert_2^2}\ 
%&\ \ \ \ \ \ \ \ \ \ \ \ \ \ \ \ \ \ \ \ \ \ \ \ \ \ 
\leq\ \sqrt{\sigma^4M+2\sigma^2\Arrowvert C\Arrowvert_{S_{\infty}}^2}\,\Arrowvert P_{U_1}-P_{U_2}\Arrowvert_{S_2},\label{eq: besser}
\end{align}
i.e. the exponentail tail bound for the differences $Z_{P_{U_1}}^{\sigma,M,,r}- Z_{P_{U_2}}^{\sigma,M,r}$ is characterized by an interplay of Hilbert-Schmidt and spectral norm, which take over the roles of the $\ell_2$- and $\ell_{\infty}$-norms of the classical Bernstein inequality from the vector case. 
%In order to derive the concentrtion inequality, 
%an 
%$S_2-S_{\infty}$-chaining device is employed. %Let $\SS_{M,r}$ denote the set of projections $P_U$ with an $r$-dimensional subspace $U\subset\mathbb{R}^M$, that is, $\SS_{M,r}$ is the Grassmann manifold, considered as a subset of $\mathbb{R}^{M\times M}$ by the isomorphism $U\mapsto P_U$. %For any $P_U\in\SS_{M,r}$, denote 
%$$
%Z_{P_U}^{M,N,r}\ :=\ \sqrt{N}\Big(\Arrowvert P_U\mathrm{E}\Arrowvert_{S_2}^2-\E\Arrowvert P_U\mathrm{E}\Arrowvert_{S_2}^2\Big).
%$$
 We prove first the case $r\leq M-r$. Note that \begin{equation}\label{eq: upper bound}
\sup_{\pi_r,\tilde{\pi}_r\in\SS_{M,r}}\Arrowvert \pi_r-\tilde{\pi}_r\Arrowvert_{S_2}^2\ \leq\ 2r\ \ \ \text{and}\ \ \ 
\sup_{\pi_r,\tilde{\pi}_r\in\SS_{M,r}}\Arrowvert \pi_r-\tilde{\pi}_r\Arrowvert_{S_{\infty}}\leq \ 2.
\end{equation}
 Since $Z_{x}^{\sigma,M,r}$ depends continuously  on $x$, 
it holds that
$$
\E\bigg(\sup_{\tilde{\pi}_r\in\SS_{M,r}}\big( Z_{\tilde{\pi}_r}^{\sigma,M,r}-Z_{{\pi}_r}^{\sigma,M,r}\big)\bigg)\ =\ \E\bigg(\sup_{\tilde{\pi}_r\in S}\big( Z_{\tilde{\pi}_r}^{\sigma,M,r}-Z_{{\pi}_r}^{\sigma,M,r}\big)\bigg)
$$
for any countable, dense  subset  $S$ of $S_{M,r}$, and by the Theorem of monotone convergence, it is sufficient to assume subsequently that $S$ is finite. 
We define now recursively an increasing family of partitions $(\AA_n)_{n\geq 0}$ of $S$ such that $\AA_0= S$, and for $n\geq 1$ and $A\in\AA_n$
\begin{align}
& (i)~~\Arrowvert \pi_r^{(1)}-{\pi}_r^{(2)}\Arrowvert_{S_2}\leq  2^{-n}\sqrt{2r}\ \ \text{and}\ \  (ii)~~\Arrowvert \pi_r^{(1)}-{\pi}_r^{(2)}\Arrowvert_{S_{\infty}} \leq 2^{-n+1}\label{eq: 1}
\end{align}
$\forall\, \pi_r^{(1)},{\pi}_r^{(2)}\in A$, with $\AA_{n+1}\subset\AA_{n}$ for all $n\geq 0$. This can be realized as follows: For $n=0$, $\AA_{0,2}:=\AA_{0,\infty}:=\{S\}$. For any totally-bounded, pseudometric space $(\mathcal{X},d)$ and any subset $E\subset\XX$, the covering number $N(E,d,\delta)$ is the smallest number of closed $d$-balls in $\XX$ of radius $\delta$ needed to cover $E$. It is proved in \cite{sza82}, see also Pajor (1998) for a different proof, that for $1\leq r\leq M-r$ and any $C'>\xi>0$
\begin{equation}
N\big(\SS_{M,r}, d_{S_2},\xi\sqrt{r}\big)\ \leq\ \Big(\frac{C'}{\xi}\Big)^{r(M-r)}\ \ \text{and}\ \ N\big(\SS_{M,r}, d_{S_{\infty}},\xi\big)\ \leq\ \Big(\frac{C'}{\xi}\Big)^{r(M-r)}
\label{eq: covering}
\end{equation}
for some universal constant $C'>0$. Hence, $S$ can be covered with at most $(C')^{r(M-r)}$ $S_2$-balls $B_{1,2},...,B_{N_{1,2},2}$ of radius $\sqrt{2r}$, and with at most $(C')^{r(M-r)}$ $S_{\infty}$-balls $B_{1,\infty},...,B_{N_{1,\infty},\infty}$ of radius ${2}$. From such finite coverings of $S_2$- and $S_{\infty}$-balls, the partitions $\AA_{1,2}$ and $\AA_{1,\infty}$ are canonically constructed by
$$
\AA_{1,j}\ :=\ \Big\{\Big(B_{k,j}\setminus\bigcup_{1\leq l< k}B_{l,j}\Big)\cap S, k=1,...,N_{1,j}\Big\},\ \  j=2,\infty.
$$ 
For $n\geq 2$ we proceed inductively using the bounds (\ref{eq: S2}) and (\ref{eq: S-infty}) for $q=\infty$ of Lemma \ref{lemma: bounds}. Indeed, each element $A\in\AA_{n-1,j}$ is element of an $S_j$-ball in $\SS_{M,r}$ of radius $2^{-n}\sqrt{2r}$ and $2^{-n+1}$, respectively, and can be partitioned as above into $(2C)^{r(M-r)}$ subsets of balls of radius $2^{-(n+1)}\sqrt{2r}$ and $2^{-n}$, respectively. By construction, the partitions $(\AA_{n,2})_{n\geq 0}$ and $(\AA_{n,\infty})_{n\geq 0}$ are nested, and $\mathrm{card}(\AA_{n,2})\leq D^{nr(M-r)}$, $\mathrm{card}(\AA_{n,\infty})\leq D^{nr(M-r)}$ for some universal constant $D>0$. Setting now $\AA_n:=\{A_2\cap A_{\infty}:\, A_2\in\AA_{2,n}\ \text{and}\ A_{\infty}\in\AA_{n,\infty}\}$ yields some partition with the above mentioned properties (\ref{eq: 1}).  Obviously,
\begin{equation}
\mathrm{card}(\AA_{n})\leq \mathrm{card}(\AA_{n,2})\mathrm{card}(\AA_{n,\infty}).\label{eq: cardinality}
\end{equation}
%At this stage, we may procced similar to the lines of the proof of Theorem 5.1, \cite{bara10}. 
For each $n\geq 1$ and $A\in\AA_n$, let $s_n(A)$ be some arbitrarily chosen rank-$r$-projection matrix in $A$. For each $s\in\SS$ and $n\geq 1$, there exists some unique $A\in\AA_n$ with $s\in A$, and we set $\Pi_n(s)=s_n(A)$. When $n=0$, define $\Pi_0(s)=\pi_r$. Now, 
\begin{align}
\E\bigg[\sup_{s\in\SS}\bigg(Z_s^{\sigma,M,r}-Z_{\pi_r}^{\sigma,M,r}\bigg)\bigg]\ &=\ \E\bigg[\sup_{s\in\SS}\sum_{n\geq 0}\Big(Z_{\Pi_{n+1}(s)}^{\sigma,M,r}-Z_{\Pi_{n}(s)}^{\sigma,M,r}\Big)\bigg]\label{eq: decomp}\\
&\leq\ \sum_{n\geq 0}\E\bigg[\sup_{s\in\SS}\Big(Z_{\Pi_{n+1}(s)}^{\sigma,M,r}-Z_{\Pi_{n}(s)}^{\sigma,M,r}\Big)\bigg].\label{eq: van}
\end{align}
Note that the decomposition $\sum_{n\geq 0}Z_{\Pi_{n+1}(s)}^{\sigma,M,r}-Z_{\Pi_{n}(s)}^{\sigma,M,r}$ of each $Z_s^{\sigma,M,r}-Z_{\pi_r}^{\sigma,M,r}$ in (\ref{eq: decomp}) is finite since $\SS$ is finite. By construction,
$$
\sup_{s\in\SS}\Arrowvert\Pi_{n}(s)- s\Arrowvert_{S_2}\ \leq\ 2^{-n}\sqrt{2r}\ \ \ \text{and}\ \ \ \sup_{s\in\SS}\Arrowvert\Pi_{n}(s)-s\Arrowvert_{S_{\infty}}\ \leq\ 2^{-n+1},
$$
hence $
\Arrowvert\Pi_{n+1}(s)- \Pi_n(s)\Arrowvert_{S_2} \leq3\cdot  2^{-n}\sqrt{2r}$ and $\Arrowvert\Pi_{n+1}(s)-\Pi_n(s)\Arrowvert_{S_{\infty}} \leq 3\cdot 2^{-n+1}$, for every $s\in\SS$. Note that $\mathrm{card}\big(\{\Pi_{n+1}(s)-\Pi_n(s): s\in\SS\}\big)\leq\mathrm{card}(A_{n+1})\mathrm{card}(\AA_n)$, $n\in\N$. Applying now Lemma A.1, \cite{vdv96}, to each expectation within  (\ref{eq: van}), yields 
\begin{align*}
\E\bigg[\sup_{s\in\SS}&\bigg(Z_s^{\sigma,M,r}-Z_{\pi_r}^{\sigma,M,r}\bigg)\bigg]\\ 
&\lesssim\ \sum_{n\geq 0}\frac{1}{2^n}\bigg(\sqrt{2r}\sqrt{\sigma^4M+2\sigma^2\Arrowvert C\Arrowvert_{S_{\infty}}^2}\sqrt{2\log\big(N_n\big)}\ +\ 2\sigma^2\log\big(N_n\big)\bigg)\\
&\lesssim\  \sigma^2rM\ +\ \sigma r\sqrt{M}\Arrowvert C\Arrowvert_{S_{\infty}},
\end{align*}
where $N_n := \mathrm{card}(\AA_{n+1})\mathrm{card}(\AA_{n}) \leq D^{(4n+2)r(M-r)}$.
 Note further that the case $r=M$ is obvious, the case $M>r>M-r$ follows by consideration of the orthogonal complements $\Arrowvert P_U(C+\mathrm{E})\Arrowvert_{S_2}^2=\Arrowvert C+\mathrm{E}\Arrowvert_{S_2}^2-\Arrowvert P_U^{\perp}(C+\mathrm{E}) \Arrowvert_{S_2}^2$.
\hfill$\square$

\section{The main result -- a general expectation bound for non-centered Gaussian matrices}\label{subsec: iteration}
%
%\subsection{The general idea for bounding the supremum}
Corollary \ref{cor: upper} (i) provides some upper bound on $\delta_{C,M,\sigma^2, r}$ which is valid for every $M\times M$-matrix $C$, and which is achieved, for instance,  for $C=0$ and $r=1$ in the small amplitude regime, and for $C=\alpha\id$ with $\arrowvert \alpha\arrowvert$ sufficiently large in the large amplitude regime, i.e.~$\arrowvert\alpha\arrowvert\gg\sigma\sqrt{M}$.  In this section, we present some new and more refined analysis for bounding $\delta_{C,M,\sigma^2,r}$ which takes advantage of some potentially  favorable structure of $C$ - resulting in the presence of the deterministic compensation term as explained in Section \ref{sec: pre}. Some protype of matrices of "high accuracy" in the large amplitude regime is discovered and analyzed.

Our approach is motivated and explained in what follows. 
The conjecture about the possibility of improvement over Corollary \ref{cor: upper} (i) for a certain type of matrices  follows from the fact that, in contrast to the situation in Section \ref{sec: weak accuracy}, the differences 
\begin{equation}\label{eq: diff}
\Arrowvert\tilde{\pi}_r(C+\mathrm{E})\Arrowvert_{S_2}^2-\Arrowvert\pi_{r}(C+\mathrm{E})\Arrowvert_{S_2}^2
\end{equation}
are usually not centered. With $\pi_r=\pi_r(C)$ maximizing the expression $\Arrowvert\tilde{\pi}_rC\Arrowvert_{S_2}^2$ over $\SS_{M,r}$ here and subsequently, the expectation of (\ref{eq: diff}) is less or equal to zero. Depending on the structure of $C$, it may be substantially smaller than zero for some appreciable amount of rank-$r$-projections. %
%\section{Non-centered Gaussian matrices}
%\subsection{Classical chaining bounds}\label{sec: non-centered}
%
%In order to motivate the subsequent analysis, the next proposition demonstrates that an expectation bound, based on chaining techniques, is sharp in the sense that it achieves the above mentioned order for centered Gaussian matrices. Here and subsequently, $\lesssim$ means less or equal up to some constant which does not depend on the variable parameters in the expression.
%
%\begin{proposition}\label{prop: step 1}
%Let $(\mathrm{E}_{ij})_{i,j=1}^M$ be a matrix of iid Gaussian entries with variance $\sigma^2$. Then for any $1\leq r<M$,
%$$
%\E\Bigg(\sup_{\tilde{\pi}_r\in\SS_{M,r}}\Arrowvert \tilde{\pi}_r\mathrm{E}\Arrowvert_{S_2}^2-\E\Arrowvert \tilde{\pi}_r\mathrm{E}\Arrowvert_{S_2}^2\Bigg)\ \lesssim\ \sigma^2rM.
%$$
%\end{proposition}
%
%The proof is a special case of the subsequent Proposition \ref{prop: rough bound} and therefore omitted. For $r=1$, the bound is of order $\sigma^2M$, which coincides with the observation (\ref{eq: intro}). 
Consequently, for any subset $\AA_C\subset\SS_{M,r}$,
\begin{align}
\sup_{\tilde{\pi}_r\in\AA_C}&\Big(\Arrowvert \tilde{\pi}_r(C+\mathrm{E})\Arrowvert_{S_2}^2-\Arrowvert \pi_r(C+\mathrm{E})\Arrowvert_{S_2}^2\Big)\label{eq: zent}\\
&\leq
\sup_{\tilde{\pi}_r\in\AA_C}\Big(\Arrowvert \tilde{\pi}_r(C+\mathrm{E})\Arrowvert_{S_2}^2-\E\Arrowvert \tilde{\pi}_r(C+\mathrm{E})\Arrowvert_{S_2}^2\nonumber\\ 
&\ \ \ \ \ \ \ \ \  \ \ \ \ \ \ \ \ \ \ -\ \Big[\Arrowvert \pi_r(C+\mathrm{E})\Arrowvert_{S_2}^2-\E\Arrowvert \pi_r(C+\mathrm{E})\Arrowvert_{S_2}^2\Big]\Big)-\Delta_{\AA_C},\nonumber
\end{align}
with
$
\displaystyle\Delta_{\AA_C}\ 
:=\ \inf_{\tilde{\pi}_r\in\AA_C}\bigg( \E\Arrowvert \pi_r(C+\mathrm{E})\Arrowvert_{S_2}^2-\E\Arrowvert \tilde{\pi}_r(C+\mathrm{E})\Arrowvert_{S_2}^2\bigg).
$

Roughly speaking, the expectation of the supremum in (\ref{eq: zent}) is small if the supremum of the corresponding centered process over $\AA_C$ is small as compared to $\Delta_{\AA_C}$. The idea for the general bound is based on decomposing the Grassmann manifold along the geometric grid of slices $\AA_{C,k}$, $k\in\N$:
\begin{align*}
&\AA_{C,k}:= \bigg\{\tilde{\pi}_r\in\SS_{M,r}:\, \frac{\Arrowvert \pi_rC\Arrowvert_{S_2}^2}{2^{k+1}} < \E\Arrowvert \pi_r(C+\mathrm{E})\Arrowvert_{S_2}^2-\E\Arrowvert \tilde{\pi}_r(C+\mathrm{E})\Arrowvert_{S_2}^2 \leq  \frac{\Arrowvert \pi_rC\Arrowvert_{S_2}^2}{2^{k}}\bigg\}.
\end{align*}
Define the random variables
$\ \ \displaystyle
Y_k\ :=\ \sup_{\tilde{\pi}_r\in\AA_{C,k}} \Big(\Arrowvert \tilde{\pi}_r(C+\mathrm{E})\Arrowvert_{S_2}^2 -\ \Arrowvert \pi_r(C+\mathrm{E})\Arrowvert_{S_2}^2\Big)
$

and
$$
Y_k^0\:=\ \sup_{\tilde{\pi}_r\in\AA_{C,k}}\hspace{-2mm}\Big(\Arrowvert \tilde{\pi}_r(C+\mathrm{E})\Arrowvert_{S_2}^2 -\E\Arrowvert \tilde{\pi}_r(C+\mathrm{E})\Arrowvert_{S_2}^2 - \Big[\Arrowvert \pi_r(C+\mathrm{E})\Arrowvert_{S_2}^2-\E\Arrowvert \pi_r(C+\mathrm{E})\Arrowvert_{S_2}^2\Big]\Big).
$$
With
$\displaystyle
\ \ \tilde{Z}\ :=\ \sup_{\substack{\tilde{\pi}_r\in\SS_{M,r}:\\ \Arrowvert\tilde{\pi}_rC\Arrowvert_{S_2}=\Arrowvert{\pi}_rC\Arrowvert_{S_2}}} \Big(\Arrowvert \tilde{\pi}_r(C+\mathrm{E})\Arrowvert_{S_2}^2 -\ \Arrowvert \pi_r(C+\mathrm{E})\Arrowvert_{S_2}^2\Big),
$

we obtain the series expansion
\begin{align}
\delta_{C,M\sigma^2,r}\ &=\ \E\sup_{\tilde{\pi}_r\in\SS_{M,r}}\Big(\Arrowvert \tilde{\pi}_r(C+\mathrm{E})\Arrowvert_{S_2}^2 -\ \Arrowvert \pi_r(C+\mathrm{E})\Arrowvert_{S_2}^2\Big)\nonumber\\
 &\leq\ \E\tilde{Z}\ +\ \sum_{k\in\mathbb{N}}\E\big(0\vee Y_k\big)\nonumber \\
&\leq\ \E\tilde{Z}\ +\ \sum_{k\in\mathbb{N}}\E\Big(0\vee \big(Y_k^0-\Delta_{\AA_{C,k}}\big)\Big).\label{eq: ztilde}
\end{align}
In view of this expansion, it is clear that good bounds can be reached in principal only if $\E Y_k^0$ is small as compared to $\Delta_{\AA_{C,k}}$. The evaluation of the expectations
$$
\E\Big(0\vee \big(Y_k^0-\Delta_{\AA_{C,k}}\big)\Big),\ \ k\in\N,
$$
is however quite a difficult task in general: It requires some suitable characterization of the subsets $\AA_{C,k}$ of the Grassmannian for general matrices $C\in\R^{M\times M}$ and tight bounds on their metric entropy. Before discussing this serious issue, we present first the main result of this section.

%It turns out that the characterization of the quality of approximation is expressed in terms of the "signal-to-noise ratio" $\tilde{C}:=C/(\sigma\sqrt{M})$ rather than $C$. Correspondingly, $\tilde{\lambda}_1\geq ...\geq\tilde{\lambda}_M$ denote the singular values of $\tilde{C}$. %The next Theorem provides a general upper bound on the approximation of the reduced-rank projection for non-centered Gaussian random matrices $X=C+\mathrm{E}$ and covers the "large amplitude" regime in particular. 

 %The subsequent claim  (i) is a special case of (ii). We decided to state it explicitly for reasons of clarity.
\begin{theorem}\label{thm: general}
Let $(\mathrm{E}_{ij})_{i,j=1}^M$ be a centered matrix of independent Gaussian entries with variance $\sigma^2$. %For any $\eta>0$, define
%$$
%\CC_{\eta}\ :=\ \big\{C\in\R^{M\times M}: \tilde{\lambda}_r^2\ \geq\ \eta(1+\tilde{\lambda}_1)\big\}.
%$$ 
%Then there exists some $\eta_0>0$, such that for any $C\in\CC_{\eta}$  with $\eta\geq\eta_0$, 
Then for any $C\in\R^{M\times M}$ with $\mathrm{rank}(C)\geq r$ and $r\leq M-r$, 
the following bound holds true:
\begin{align*}
\delta_{C,M,\sigma^2,r}\ \lesssim \ \sigma^2rM\Big( \mathrm{I}+\min\big(\mathrm{II}, \mathrm{III}\big)\Big),
\end{align*}
where, with $\lambda_1\geq\lambda_2\geq ...\geq\lambda_M$ denoting the singular values of $C$,
\begin{align*}
\mathrm{I}\ &=\  \min\bigg(\frac{\lambda_1^2}{{\lambda}_r^2}, 1+\frac{\lambda_1}{\sigma\sqrt{M}}\bigg), \\
\mathrm{II}\ &=\ \bigg(\frac{\frac{1}{r}\sum_{i=r+1}^{2r}\lambda_i^2}{\lambda_r^2}\bigg)^{1/2}\cdot\frac{{\lambda}_1}{\sigma\sqrt{M}},\ \ \ \text{and}\\
\mathrm{III}\ &=\ \frac{\lambda_1^2}{\lambda_r^2-\lambda_{r+1}^2}\ \ \text{if}\ \lambda_{r+1}<\lambda_r,\ \ \text{and}\ \ \mathrm{III}=\infty\ \ \text{else}.
\end{align*}
\end{theorem}

In particular for the case $r=1$, the bound applies to any $0\not =C\in\R^{M\times M}$. Some immediate consequence of Theorem \ref{thm: general} is the following.
%\begin{corollary}[Prototype of high accuracy in the large amplitude regime]Let $C_{\alpha}\in\R^{M\times M}$ and let $U\Lambda_{\alpha}V'$ be some singular value decomposition of $C_{\alpha}$, where $\Lambda_{\alpha}=\alpha\id_r$. Then 
%$$
%\delta_{C_{\alpha},M,\sigma^2,r}\ \lesssim\ \sigma^2rM
%$$
%with some fucntion $\gamma(\cdot)\leq \log(1+r)$ and $\lim_{\alpha\rightarrow\infty}\gamma(M,r,\alpha)=1$.
%\end{corollary}
\begin{corollary}\label{cor: step}Let $C_{\alpha,\beta,r}\in\R^{M\times M}$ with singular values $\lambda_i=\alpha$ for $i\leq r$ and\linebreak $\lambda_i=\beta\leq\alpha$ for $i > r$. As usual, set $c/0:=\infty$ for any $c>0$. Then 
$$
\delta_{C_{\alpha,\beta,r},M,\sigma^2,r}\ \lesssim\ \sigma^2rM\bigg\{1\ +\  \min\bigg(\frac{\alpha^2}{\alpha^2-\beta^2},\,\frac{\beta}{\sigma\sqrt{M}}\bigg)\bigg\}\ \ \text{for \underline{every}}\ \ 0\leq \beta\leq \alpha.
$$
\end{corollary}
\begin{remark} Corollary (\ref{cor: step}) covers the two extreme cases:
\begin{itemize}
\item[(i)] Prototype of high accuracy in the "large amplitude"-regime  $\alpha\gg \sigma\sqrt{M}$:
$$
\delta_{C_{\alpha,0,r},M,\sigma^2,r}\ \lesssim\ \sigma^2rM\ \ \ \text{for \underline{every}}\ \ \alpha\geq 0.
$$
\item[(ii)] Prototype of weak accuracy in the "large amplitude"-regime $\alpha\gg \sigma\sqrt{M}$:
the bound of Corrollary \ref{cor: step} approaches the unimprovable upper bound for matrices $C$ of constant singular value spectrum as $\beta\rightarrow\alpha$, see Proposition \ref{prop: stability}. 
\end{itemize}
Note that the upper bound in case $C=0$ is covered by Proposition \ref{prop: rough bound}.
\end{remark}

%\begin{remark}(Spectral gap in $r$) The condition of a spectral gap in $r$ is equivalent to uniqueness of $\pi_r$ and is assumed for technical reasons, as it allows for some very clean description  of the slices $A_{C,k}$ in Subsection \ref{subsec: slices}. Extensions are straightforward but very technical. The bound is getting larger the more massive the set $$
%\underset{~~\tilde{\pi}_r\in\SS_{M,r}}{\mathrm{Argmax}}\,\Arrowvert\tilde{\pi}_rC\Arrowvert_{S_2}^2$$  
%is, which also determines the value of $\E\tilde{Z}$ in (\ref{eq: ztilde}). We discuss in particular the non-uniqueness of $\pi_r$ at the end of the next subsection.
%\end{remark}

The proof of Theorem \ref{thm: general} is deferred to Section \ref{sec: theorem main}. Subsection \ref{subsec: slices} deals with the description of the sets $\GG_{M,r}(\delta):=\{\tilde{\pi}_r\in\SS_{M,r}: \Arrowvert \pi_rC\Arrowvert_{S_2}^2 -\Arrowvert\tilde{\pi}_r C\Arrowvert_{S_2}^2\leq\delta\}$, which characterize the slices $\AA_{C,k}= \GG_{M,r}\big(2^{-k}\Arrowvert\pi_r C\Arrowvert_{S_2}^2\big)\setminus \GG_{M,r}\big(2^{-{(k+1)}}\Arrowvert\pi_r C\Arrowvert_{S_2}^2\big)$. It is shown that  these sets can be approximated by sets of  very simple geometric structure. In case of some substantial spectral gap in $r$, this approximation is very tight. Sharp bounds on their metric entropy are derived in Subsection \ref{sec: rand}. The final arguments differ slightly from from the description at the beginning of this Section. They are given in Subsection \ref{subsec: general}.

\section{Lower bounds}\label{sec: lower bounds}
The best possible upper bound for the accuracy of aproximation in Theorem \ref{thm: general} is of the order $\sigma^2rM$, which is attained for matrices $C=0$ or $C=\alpha\id_r$ in the large amplitude regime. The question arises whether this bound is sharp, i.e.~whether it  indicates some fundamental limit on the accuracy of approximation. For fixed $A\in\R^{M\times M}$ with $\Arrowvert A\Arrowvert_{S_{\infty}}=1$ and some arbitrary real number $\alpha\in\R$,
inspection of $\delta_{\alpha A,M,\sigma^2,r}$ shows
\begin{align}
\E\big\Arrowvert \hat{\pi}_r& (\alpha A+\mathrm{E})\big\Arrowvert_{S_2}^2 \ -\ \E\big\Arrowvert \pi_r (\alpha A+\mathrm{E})\big\Arrowvert_{S_2}^2\nonumber\\ 
&=\ \E\bigg(\sup_{\tilde{\pi}_r\in\SS_{M,r}}\big\Arrowvert\tilde{\pi}_r\mathrm{E}\big\Arrowvert_{S_2}^2-\big\Arrowvert{\pi}_r\mathrm{E}\big\Arrowvert_{S_2}^2 +2\alpha\,\tr\Big(\mathrm{E}'(\tilde{\pi}_r-\pi_r)A\Big)\label{eq: alpha-1}\\
&\ \ \ \ \ \ \ \ \ \ \ \ \ \ \ \ \ \ \ \ \ -\alpha^2\Big(\big\Arrowvert{\pi}_rA\big\Arrowvert_{S_2}^2-\big\Arrowvert \tilde{\pi}_rA\big\Arrowvert_{S_2}^2\Big)\bigg).\label{eq: alpha-2}
\end{align}
Now, $\Arrowvert{\pi}_rC\Arrowvert_{S_2}^2-\Arrowvert \tilde{\pi}_rC\Arrowvert_{S_2}^2 \geq 0$ since $\pi_r$ optimizes the Hilbert-Schmidt norm, and the dependence on $\alpha$ in the deterministic compensation term in line (\ref{eq: alpha-2}) is quadratic while it is only linear in the stochastic part in line (\ref{eq: alpha-1}). So, for any fixed $\sigma^2$, $r$, $M$ and $A$, one may wonder whether the accuracy  of approximation $\E\Arrowvert \hat{\pi}_r (\alpha \id_r+\mathrm{E})\Arrowvert_{S_2}^2  - \E\Arrowvert \pi_r (\alpha \id_r+\mathrm{E})\Arrowvert_{S_2}^2$ tends even to zero as $\arrowvert \alpha\arrowvert$ goes to infinity  if $A$ is suitably chosen, like, for instance, $A=\id_r$.
In this section, we demonstrate that this is not the case. We provide some complete proof of the conjecture 
\begin{equation}\label{eq: conjecture}
\inf_{C\in\R^{M\times M}}\delta_{C,M,\sigma^2,r}\ \gtrsim\ \delta_{0,M,\sigma^2,r}
\end{equation}
in case $r=1$. For  $r>1$ we present some partial solution in the large amplitude regime.

\subsection{The universal lower bound for $r=1$}
\begin{theorem}\label{thm: r=1}
Let $(\mathrm{E}_{ij})_{i,j=1}^M$ be a centered matrix of independent Gaussian entries with variance $\sigma^2$. Then there exists some $M_0\in\N$, independent of $\sigma^2$, such that
\begin{equation}
\inf_{C\in\R^{M\times M}}\delta_{C,M,\sigma^2,1}\ \gtrsim\ \delta_{0,M,\sigma^2,1}
\end{equation}
for any $M\geq M_0$.
\end{theorem}

\paragraph{\sc Proof}In view of Lemma \ref{prop: stability2} and its subsequent remark, it is sufficient to prove that for every $\beta>0$, there exists some constant $c_{\beta}>0$, independent of $\sigma^2$ and $M$, such that 
\begin{align}\label{eq: lower sharp}
\inf_{\substack{C\in\R^{M\times M}:\\ \Arrowvert C\Arrowvert_{S_{\infty}}\geq \beta\sigma\sqrt{M}}}\  \geq\ c_{\beta}\sigma^2(M-1)\ \ \text{for all}\ \  M\geq 2.
\end{align}
\noindent
Let $C=\sum_{i=1}^M\lambda_iU_iV_i'$ denote some singular value decomposition of $C$, and define $\pi_1:=U_1U_1'$. 
Since for any $1\leq s<M$ and $\pi_s,\tilde{\pi}_s\in\SS_{M,s}$, $\pi_s-\tilde{\pi}_s=(\id-\pi_s)\tilde{\pi}_s-\pi_s(\id -\tilde{\pi}_s)$ is an orthogonal decomposition and $\Arrowvert \pi_s(\id -\tilde{\pi}_s)\Arrowvert_{S_2}^2=\Arrowvert (\id -\pi_s)\tilde{\pi}_s\Arrowvert_{S_2}^2$, observe that
\begin{align*}
\Arrowvert\pi_1 C\Arrowvert_{S_2}^2-\Arrowvert \tilde{\pi}_1C\Arrowvert_{S_2}^2\ &\leq\ \Arrowvert\pi_1 C\Arrowvert_{S_2}^2-\Arrowvert \tilde{\pi}_1\pi_1C\Arrowvert_{S_2}^2\\
&=\ \lambda_1^2\Big(\Arrowvert\pi_1\Arrowvert_{S_2}^2-\Arrowvert \tilde{\pi}_1\pi_1\Arrowvert_{S_2}^2\Big)\\
&=\ \lambda_1^2\Arrowvert (\id -\tilde{\pi}_1)\pi_1\Arrowvert_{S_2}^2\\
&=\  \lambda_1^2\Arrowvert \tilde{\pi}_1(\id-\pi_1)\Arrowvert_{S_2}^2 =\ \frac{\lambda_1^2}{2}\Arrowvert \tilde{\pi}_1-\pi_1\Arrowvert_{S_2}^2.
\end{align*}
Consequently,
\begin{align}
\delta_{C,M,\sigma^2,1}\ &=\ \E\big\Arrowvert \hat{\pi}_1 (C+\mathrm{E})\big\Arrowvert_{S_2}^2 \ -\ \E\big\Arrowvert \pi_1 ( C+\mathrm{E})\big\Arrowvert_{S_2}^2\nonumber\\ 
&=\ \E\bigg(\sup_{\tilde{\pi}_1\in\SS_{M,r}}\bigg\{\big\Arrowvert\tilde{\pi}_1\mathrm{E}\big\Arrowvert_{S_2}^2-\big\Arrowvert{\pi}_1\mathrm{E}\big\Arrowvert_{S_2}^2 +2\,\tr\Big(\mathrm{E}'(\tilde{\pi}_1-\pi_1)C\Big)\nonumber\\
&\ \ \ \ \ \ \ \ \ \ \ \ \ \ \ \ \ \ \ \ \ \ \ \ \ \ \ \ \ \ \ \ \ -\Big(\big\Arrowvert{\pi}_1C\big\Arrowvert_{S_2}^2-\big\Arrowvert \tilde{\pi}_1C\big\Arrowvert_{S_2}^2\Big)\bigg\}\bigg)\nonumber\\
&\geq\ \E\bigg(\sup_{\tilde{\pi}_1\in\SS_{M,r}}\bigg\{\big\Arrowvert\tilde{\pi}_1\mathrm{E}\big\Arrowvert_{S_2}^2-\big\Arrowvert{\pi}_1\mathrm{E}\big\Arrowvert_{S_2}^2 +2\,\tr\Big(\mathrm{E}'(\tilde{\pi}_1-\pi_1)C\Big)\nonumber\\
&\ \ \ \ \ \ \ \ \ \ \ \ \ \ \ \ \ \ \ \ \ \ \ \ \ \ \ \ \ \ \ \ \ -\frac{\lambda_1^2}{2}\Arrowvert\pi_1-\tilde{\pi}_1\Arrowvert_{S_2}^2\bigg\}\bigg)\nonumber\\
&\geq\  \E\bigg(\sup_{\tilde{\pi}_1\in\SS_{M,1}(\delta)}\bigg\{\big\Arrowvert\tilde{\pi}_1\mathrm{E}\big\Arrowvert_{S_2}^2-\big\Arrowvert{\pi}_1\mathrm{E}\big\Arrowvert_{S_2}^2 +2\,\tr\Big(\mathrm{E}'(\tilde{\pi}_1-\pi_1)C\Big)\label{eq: alpha-1b}\\
&\ \ \ \ \ \ \ \ \ \ \ \ \ \ \ \ \ \ \ \ \ \ \ \ \ \ \ \ \ \ \ \ \ -\frac{\lambda_1^2}{2}\Arrowvert\pi_1-\tilde{\pi}_1\Arrowvert_{S_2}^2\bigg\}\bigg)\label{eq: alpha-2b} 
\end{align}
for any subset
$$
\SS_{M,1}(\delta)\ :=\ \Big\{\tilde{\pi}_1\in\SS_{M,1}:\, \Arrowvert\tilde{\pi}_1-\pi_1\Arrowvert_{S_2}\leq\sqrt{2\delta}\Big\},\ \ \ \delta>0.
$$
The idea of the proof is to choose $\delta=\delta(M,\sigma^2,\lambda_1,\beta)$ in some specific way in order to guarantee that the deterministic compensation term in (\ref{eq: alpha-2b}) is lower bounded by $-\lambda_1^2\delta$,  and to pick afterwards some suitable projection in dependence of $C$ and $\mathrm{E}$ out of this class which realizes the lower bound in expectation. Since $\check{\pi}_{r,U}=U'\tilde{\pi}_rU\in\SS_{M,r}$ for any $\tilde{\pi}_r\in\SS_{M,r}$ and $U\in O(M)$, $\check{\mathrm{E}}_{U,V}=U'\mathrm{E} V=_{\DD}\mathrm{E}$ for any two fixed orthogonal matrices $U$ and $V$, and $\Arrowvert U'AV\Arrowvert_{S_2}^2=\Arrowvert A\Arrowvert_{S_2}^2$ for any $A\in\R^{M\times M}$, i.e., $$
\big\Arrowvert \tilde{\pi}_r\big(U\Lambda V'+\mathrm{E}\big)\big\Arrowvert_{S_2}^2\ =\ \big\Arrowvert \check{\pi}_{r,U}\big(\Lambda +\check{\mathrm{E}}_{U,V}\big)\big\Arrowvert_{S_2}^2,
$$
we may and do assume subsequently without loss of generality that $U=V=\id$ and  $\pi_1=\id_1$ in particular.
With $e_1,e_2,...,e_M$ denoting the canonical basis vectors of $\R^M$, every projection matrix in $\SS_{M,1}$ can be written as
$$
\Big(\sum_{i=1}^M\gamma_ie_i\Big)\Big(\sum_{i=1}^M\gamma_i e_i'\Big)\ =:\ \tilde{\pi}_{1,\gamma},\ \ \text{with}\ \sum_{i=1}^M\gamma_i^2=1.
$$
In order bound $\delta_{C,M,\sigma^2,1}$ from below by means of (\ref{eq: alpha-1b}) -- (\ref{eq: alpha-2b}), define
$$
\delta_{**}\ :=\ d\frac{(M-1)\sigma^2}{2\lambda_1^2}
$$
for some constant $d\in (0,\min(\beta^2,1))$ to be chosen later. 
Note that $d$ shall be chosen independently of $M$, $\sigma^2$, $\lambda_1$ and $C$, but is allowed to depend on $\beta$ only. Furthermore, $\delta^*\leq 1/2$ because $\lambda_1\geq\beta\sigma\sqrt{M}$. Now, since
$$
\tilde{\pi}_{1,\gamma}\id_1\ =\ \left(\begin{array}{cccc}
\gamma_1^2 & 0&\cdots& 0\\
\gamma_1\gamma_2 & 0 & \cdots & 0\\
\vdots & 0 & \ddots & \vdots\\
\gamma_1\gamma_M & 0  &\cdots &0 
                                              \end{array}\right),
$$
the constraint 
$$\tilde{\pi}_{1,\gamma}\in \SS_{M,1}(\delta_{**})\ \Leftrightarrow\ \Arrowvert \tilde{\pi}_{1,\gamma}-\id_1\Arrowvert_{S_2}^2\leq 2\delta_{**}\ \Leftrightarrow\ \Arrowvert (\id-\tilde{\pi}_{1,\gamma})\id_1\Arrowvert_{S_2}^2\leq\delta_{**}$$ 
translates into 
$$
\Big((1-\gamma_{1}^2)^2+\sum_{i=2}^M\gamma_1^2\gamma_i^2\ =\ (1-\gamma_1^2)^2+\gamma_1^2(1-\gamma_1^2)\ =\Big)\ \ 1-\gamma_1^2\ \leq\  \delta_{**}.
$$
With the choice 
$$
\gamma_1^*\ :=\ \sqrt{1-\delta_{**}}\ \ \  \text{and}\ \ \  \gamma_i^*\ :=\ \mathrm{sign}(\mathrm{E}_{i1})\,\delta_{**}^{1/2}/\sqrt{M-1}\ \text{ for $i=2,...,M$},
$$ 
it holds that $\Arrowvert \tilde{\pi}_{1,\gamma^*}-\id_1\Arrowvert_{S_2}^2\leq 2\delta_{**}$, i.e.~$\tilde{\pi}_{1,\gamma^*}$ belongs to $\SS_{M,1}(\delta_{**})$. Together with $\delta=\delta_{**}$ in (\ref{eq: alpha-1b}) -- (\ref{eq: alpha-2b}), this yields
\begin{equation}\label{eq: direct}
\delta_{C,M,\sigma^2,1}\ \geq\ \E\bigg(\big\Arrowvert\tilde{\pi}_{1,\gamma^*}\mathrm{E}\big\Arrowvert_{S_2}^2-\big\Arrowvert{\id}_1\mathrm{E}\big\Arrowvert_{S_2}^2 +2\,\tr\Big(\mathrm{E}'(\tilde{\pi}_{1,\gamma*}-\id_1)C\Big)\ -\ d\sigma^2(M-r)\bigg).
\end{equation}
We evaluate the expressions within the expectation separately.
For any $1\leq i,j\leq M$,
\begin{align*}
\big(\tilde{\pi}_{1,\gamma}\mathrm{E}\big)_{ij}\ =\ \sum_{l=1}^M\gamma_i\gamma_l\mathrm{E}_{lj}.
\end{align*}
Since $\sum_{i=1}^M{\gamma_i^*}^2=1$, $\mathrm{sign}(\mathrm{E}_{i1})$ and $\arrowvert \mathrm{E}_{i1}\arrowvert$ are stochastically independent, and $\mathrm{E}_{ij}$, $1\leq i, j\leq M$, are independent by assumption with $\E(\mathrm{E}_{ij})=0$ and $\E\,\mathrm{sign}(\mathrm{E}_{ij}) =0$,
\begin{align*}
\E\Arrowvert \tilde{\pi}_{1,\gamma^*}\mathrm{E}\Arrowvert_{S_2}^2\ &=\ \sum_{i,j=1}^M\E\bigg(\sum_{l=1}^M\gamma^*_i\gamma^*_l\mathrm{E}_{lj}\bigg)^2\\ 
&=\ \sum_{i,j=1}^M\sum_{l=1}^M{\gamma^*_i}^2{\gamma^*_l}^2\E\big(\mathrm{E}_{lj}^2\big)\ +\ \sum_{\substack{l,l'\geq 2:\\l\not= l'}}^M \arrowvert\gamma_l^*\gamma_{l'}^*\arrowvert\E\big\arrowvert\mathrm{E}_{l1}\mathrm{E}_{l'1}\big\arrowvert \\
&=\ M\sigma^2\ +\ \sum_{\substack{l,l'\geq 2:\\l\not= l'}}^M \arrowvert\gamma_l^*\gamma_{l'}^*\arrowvert\E\big\arrowvert\mathrm{E}_{l1}\mathrm{E}_{l'1}\big\arrowvert.
\end{align*}
Therefore,
\begin{equation}\label{eq: quad 0}
\E\Big(\big\Arrowvert\tilde{\pi}_{1,\gamma^*}\mathrm{E}\big\Arrowvert_{S_2}^2-\big\Arrowvert{\id}_1\mathrm{E}\big\Arrowvert_{S_2}^2\Big)\ \geq\ 0.
\end{equation}
Next, we decompose
\begin{equation}
2\,\tr\Big(\mathrm{E}'\tilde{\pi}_{1,\gamma*}C\Big)\ =\ 2\,\tr\Big(\mathrm{E}'\tilde{\pi}_{1,\gamma*}\id_1C\Big)\ +\ 2\,\tr\Big(\mathrm{E}'\tilde{\pi}_{1,\gamma*}(\id-\id_1)C\Big).
\end{equation}
In order to check $\E\,\tr\big(\mathrm{E}'\tilde{\pi}_{1,\gamma^*}(\id-\id_1)C\big)=0$, it is sufficient to notice that all entries of the matrix $\mathrm{E}'\tilde{\pi}_{1,\gamma^*}(\id-\id_1)$ have expectation $0$. Indeed, its first column is equal to zero, and one easily verifies  for the remaining indices $1\leq i\leq M$, $2\leq j\leq M$ that $\E\big(\mathrm{E}_{li}\gamma^*_j\gamma^*_l\big)=0$ for every $1\leq l\leq M$, hence
$$
\Big(\E\big(\mathrm{E}'\tilde{\pi}_{1,\gamma^*}(\id-\id_1)\big)\Big)_{ij}\ =\ \sum_{l=1}^M\E\big(\mathrm{E}_{li}\gamma^*_j\gamma^*_l\big)\ =\ 0.
$$
Together with (\ref{eq: quad 0}) and $\id_1C=\lambda_1\id_1$, (\ref{eq: direct}) reduces to
\begin{align*}
\delta_{C,M,\sigma^2,1}\ &\geq\ \E\bigg(2\,\tr\Big(\mathrm{E}'(\tilde{\pi}_{1,\gamma*}-\id_1)\id_1C\Big)\ -\ d\sigma^2(M-1)\bigg)\\
&=\ \E\bigg(2\lambda_1\,\tr\Big(\mathrm{E}'\tilde{\pi}_{1,\gamma*}\id_1\Big)\ -\ d\sigma^2(M-1)\bigg)\\
&=\ \sigma\cdot2\lambda_1\,\E\left({\gamma_1^*}^2(\mathrm{E}_{11}/\sigma)+\sum_{i=2}^M\gamma_1^*\gamma_i^*(\mathrm{E}_{i1}/\sigma)\right)\ -\ d\sigma^2(M-1)\\
&=\ (1-\delta_{**})^{1/2}\delta_{**}^{1/2}(M-1)^{-1/2}\sigma2\lambda_1\sum_{i=2}^M\E\arrowvert\mathrm{E}_{i1}/\sigma\arrowvert\ -\ d\sigma^2(M-1)\\ 
&\geq \ \sqrt{d}\,\sigma^2(M-1)\frac{2}{\sqrt{2\pi}}\ -\ d\sigma^2(M-1).
\end{align*}
Choosing now some $d\in \big(0,\min(\beta^2,1)\big)$ such that $2\sqrt{d}/\sqrt{2\pi}-d>0$ proves (\ref{eq: lower sharp}).
\hfill$\square$

\begin{remark}
Theorem \ref{thm: r=1} answers to  question (A) from the Introduction in the negative.
\end{remark}

\subsection{A partial solution to the conjecture (\ref{eq: conjecture}) for $r>1$ in the large amplitude regime}
The specific construction of the (random) projection $\tilde{\pi}_{1,\gamma^*}$ in the proof from the previous paragraph cannot canonically be extended to arbitrary $r>1$. For the result in this subsection, we use finally a different approach based on abstract lower bounds on suprema of Gaussian processes which applies to any $r\leq M-r$. For the prototype $C_{\alpha}=\alpha\id_r$ of high accuracy in the large amplitude regime $\arrowvert \alpha\arrowvert \gg\sigma{\sqrt{M}}$ of Theorem \ref{thm: general}, the following result is deduced. Its substantially more involved extension to a non-asymptotic optimal lower bound for general matrices $C$ may also involve the Sudakov-type minoration for Gaussian chaos processes (\cite{tala92}).

\begin{theorem}\label{thm: lower bounds}
Let $(\mathrm{E}_{ij})_{i,j=1}^M$ be a centered matrix of independent Gaussian entries with variance $\sigma^2$. Let $C_{\alpha,s}\in\R^{M\times M}$ with singular value decomposition $U\Lambda_{\alpha,s}V'$, where $\Lambda_{\alpha,s}=\alpha\id_s$ with $0<\alpha\in\R$ and $1\leq s<M$. Then
\begin{equation}\label{eq: asymptotic LB}
\underset{{\alpha}\rightarrow\infty}{\lim\inf}\,\max_{s\in\{r,M-r\}}\delta_{C_{\alpha,s},M,\sigma^2,s}\ \gtrsim\ \sigma^2r(M-r).
\end{equation}
\end{theorem}

The proof of Theorem \ref{thm: lower bounds} is deferred to Section \ref{sec: lower bounds 3}.

\begin{remark}
In view of the polar decomposition of $(\id-\tilde{\pi}_r)\pi_r$ and $\tilde{\pi}_r(\id-\pi_r)$ which shows in particular that these two matrices have the same singular values, we conjecture that the bound holds true also without the maximum over $s\in\{r,M-r\}$, but do not have a rigorous and elegant proof for it yet. Note that this maximum is redundant if $r=M/2$, $M\in 2\N$.
\end{remark}

\section{Consequences on statistical estimation problems}\label{sec: discussion} Let $X=C+\mathrm{E}$ as described in (\ref{eq: model 1}). Let $\hat{\lambda}_1\geq\hat{\lambda}_2\geq ...\geq \hat{\lambda}_M$ and  ${\lambda_1}\geq{\lambda_2}\geq ...\geq {\lambda}_M$ denote the singular values of $X$ and $C$, respectively. Recall that $\sum_{i=1}^r\hat{\lambda}_i^2=\Arrowvert\hat{\pi}_rX\Arrowvert_{S_2}^2$ and $\sum_{i=1}^r\lambda_i^2=\Arrowvert \pi_rC\Arrowvert_{S_2}^2$ with the rank-$r$-projections $\hat{\pi}_r$ and $\pi_r$ as defined in the Introduction.

\subsection{The largest singular value}\label{subsec: largest}
We begin with the simplest example of estimating  $\lambda_1^2$, the largest eigenvalue of $C'C$, based on the observation $X=C+\mathrm{E}$. As explained in the introduction, the maximal eigenvalue of $X'X$ is positively biased as an estimator for $\lambda_1^2$, because
$$
\E\,\hat{\lambda}_1^2\ =\ \E\Arrowvert\hat{\pi}_1X\Arrowvert_{S_2}^2\ \geq\ \E\Arrowvert{\pi}_1X\Arrowvert_{S_2}^2\ =\ \lambda_1^2\ +\ \sigma^2M.
$$
Therefore, one immediate improvement is to consider $\hat{s}:=\hat{\lambda}_1^2-\sigma^2M$ as an estimator for $\lambda_1^2$.  As Theorem \ref{thm: r=1} reveals for the particular case $r=1$,  
\begin{equation}\label{eq: beispiel}
\E\,\hat{s}-\lambda_1^2\ =\ \E\,\hat{\lambda}_{1}^2-\sigma^2 M -\lambda_1^2
\end{equation}
is stricly positive and bounded away from zero, uniformly over $C\in\R^{M\times M}$. As a consequence of Corollary \ref{cor: upper} and Theorem \ref{thm: r=1},  $$\delta_{C,M,\sigma^2,1}\ =\ \E\,\hat{s}-\lambda_1^2
\ \in\ \Big[c_1\sigma^2M,\,c_2\big(\sigma^2M+\sigma \sqrt{M}\Arrowvert C\Arrowvert_{S_{\infty}}\big)\Big]
$$%    ranges between $c_1\sigma^2M$ and $c_2\big(\sigma^2M+\sigma \sqrt{M}\Arrowvert C\Arrowvert_{S_{\infty}}\big)$  
for some universal real constants $c_1,c_2>0$ which do not depend on $M,\sigma^2$ and $C$, and it follows from (\ref{eq: intro 2}) for $C=0$ and Proposition \ref{prop: stability} (ii) for $C=\alpha\id$ that these bounds cannot be improved in general. Hence, one message of our analysis is: 
\begin{itemize}\item[] The quantity $\sigma^2M$ always underestimates the bias $\E\big(\hat{\lambda}_1^2-\lambda_1^2\big)$ by at least some universal factor strictly larger than $1$, independently on how favorable the matrix $C$ is.
 \end{itemize}
In other words, even after correction by $\sigma^2M$, the difference \begin{equation}\label{eq: difference 6}\E\,\hat{\lambda}_1^2-\lambda_1^2-\sigma^2M\end{equation} remains of the same order $\sigma^2M$ at least, independently of $C$. Moreover, there exist matrices $C$ for which (\ref{eq: beispiel}) is not smaller in order than $\sigma^2M+\sigma \sqrt{M}\Arrowvert C\Arrowvert_{S_{\infty}}$. That is, large amplitude $\Arrowvert C\Arrowvert_{S_{\infty}}$ never improves (in order) the acuracy  as compared to $C=0$, but it may result in substantially worse accuracy of approximation.  
Therefore, some further consequence is that small magnitude of $\sigma^2M$ is necessary but far from being sufficient for the bias of $\hat{s}$ in (\ref{eq: beispiel}) to be small. The worst case error is non-asymptotically quantified in terms of $\Arrowvert C\Arrowvert_{S_{\infty}}, \sigma^2$ and $M$ in Corollary \ref{cor: upper} (i). Theorem \ref{thm: general} describes more precisely the effect of the shape of the singular value spectrum on the accuracy of approximation. For example, if $C=\alpha\id$, then (\ref{eq: difference 6})
grows like $\arrowvert\alpha\arrowvert\sigma \sqrt{M}$ as
$\arrowvert\alpha\arrowvert\rightarrow\infty$, cf.~Proposition \ref{prop: stability} (ii), for any fixed $\sigma^2$ and $M$. On the other hand, if $C=\alpha\id_1$, then (\ref{eq: difference 6}) remains bounded by some universal constant times $\sigma^2M$, independently of $\alpha$.  The same holds true for the full rank matrix $C=\alpha\id_1+\alpha\id$. Consequently, Theorem \ref{thm: general} demonstrates that large amplitude $\Arrowvert C\Arrowvert_{S_{\infty}}$ does not necessarily result in worse accuracy of approximation as compared to the case $C=0$, and discovers some  prototypes of high accuracy in the large amplitude regime. 
Similar conclusions  for $r>1$, i.e.~statistics of the form $\sum_{i=1}^r\hat{\lambda}_i^2$, are valid as well.

\subsection{Relative quantities} This subsection is devoted to the consequences of our results on relative quantities as described in the introduction. Consider, for instance, the ratio
$$
{t}_r\ :=\ \frac{\Arrowvert\pi_rC\Arrowvert_{S_2}^2}{\Arrowvert C\Arrowvert_{S_2}^2}\ =\ \frac{\sum_{i=1}^r\lambda_i^2}{\sum_{i=1}^M\lambda_i^2}.
$$ 
Assuming $\Arrowvert C\Arrowvert_{S_2}^2$ to be known, some natural candidate estimator of $t_r$ is
$$
\hat{t}_{r}\ :=\ \frac{\Arrowvert\hat{\pi}_rX\Arrowvert_{S_2}^2-\sigma^2rM}{\Arrowvert C\Arrowvert_{S_2}^2}.
$$
In this situation, increasing amplitude of $C$ always results in smaller bias of the estimator $\hat{t}_r$. Suppose that $C=\alpha A$ for some real number $\alpha\in\R$ and an $M\times M$-matrix $A$ with $\Arrowvert A\Arrowvert_{S_{\infty}}=1$. Note that $\Arrowvert A\Arrowvert_{S_2}^2\geq \Arrowvert A\Arrowvert_{S_{\infty}}^2=1$. Then Corollary \ref{cor: upper} (i) yields
$$
\sup_{\substack{A\in\R^{M\times M}:\\
\mathrm{rank}(A)\geq r,\ \Arrowvert A\Arrowvert_{S_{\infty}}=1}}\E\big(\hat{t}_{r} - t_r\big)\ \lesssim\ \frac{\sigma^2rM}{\alpha^2}\bigg\{1+\frac{\arrowvert \alpha\arrowvert }{\sigma\sqrt{M}}\bigg\}.
$$
 For every fixed $A$, $\sigma^2$, $r$ and $M$, the bias $\E(\hat{t}_{r} - t_r)$ tends to zero as $\arrowvert\alpha\arrowvert\rightarrow\infty$, in contrast to the situation for the absolute difference in the previous Subsection \ref{subsec: largest}. Note that it decreases of the order $\arrowvert\alpha\arrowvert^{-1}$ at least and of the order $\alpha^{-2}$ at most, depending on the shape of the singular value spectrum of $A$, cf.~Theorem \ref{thm: general}. So, whereas, for every fixed $A$, $\sigma^2$, $r$  and $M$, the {\it absolute difference} as described in Subsection \ref{subsec: largest} for $r=1$ cannot get closer to zero as the amplitude $\arrowvert\alpha\arrowvert$ of $C=\alpha\,A$ increases, independently on how favorable the matrix $C$ may be, the {\it relative difference} $\E(\hat{t}_r-t_r)$ tends to zero as the amplitude goes to infinity for every matrix $C=\alpha\,A$. The shape of the singular value spectrum however clearly influences the accuracy of approximation, in the same fashion as explained in Subsection \ref{subsec: largest}, cf. Theorem \ref{thm: general}.% the order of which in terms of the amplitude $\alpha$  ranges between  $\arrowvert\alpha\arrowvert^{-1}$ at least and of the order $\alpha^{-2}$ at most.

\subsection{Quadratic functionals of low-rank matrices}
One natural candidate for estimating $\Arrowvert C\Arrowvert_{S_2}^2$, based on the observation $X=C+\mathrm{E}$ described by (\ref{eq: model 1}), is the unbiased estimator
$
\Arrowvert X\Arrowvert_{S_2}^2 - \sigma^2M^2
$.
Simple calculation yields 
\begin{equation}\label{eq: variance 1}
\mathrm{Var}\Big(\Arrowvert X\Arrowvert_{S_2}^2-\sigma^2M^2\Big)\ =\ 2\sigma^4M^2+4\sigma^2\Arrowvert C\Arrowvert_{S_2}^2.
\end{equation}
One disadvantage of this estimator is its large variance for large values of $M$: it depends quadratically on the dimension. %Some dimension reduction assumption which found recently increasing interest in the statistical literature  is the "low rank hypothesis" 
%$$
%\mathrm{rank}(C)\ll M.
%$$
If $r=\mathrm{rank}(C)<M$, then the matrix $C$ can be fully characterized by $(2M-r)r$ parameters as can be seen by the singular value decomposition. That is, if $r\ll M$, the intrinsic dimension of the problem is of the order $rM$ rather than $M^2$. Now observe that for every matrix $C$ with $\mathrm{rank}(C)=r$,
$$
\Arrowvert C\Arrowvert_{S_2}^2\ =\ \Arrowvert \pi_rC\Arrowvert_{S_2}^2.
$$
Elementary calculation reveals that $\Arrowvert\pi_rX\Arrowvert_{S_2}^2-\sigma^2rM$ unbiasedly estimates $\Arrowvert C\Arrowvert_{S_2}^2$, and
\begin{equation}\label{eq: variance 2}
 \mathrm{Var}\Big(\Arrowvert\pi_rX\Arrowvert_{S_2}^2-\sigma^2rM\Big)\ =\ 2\sigma^4rM\, +\, 4\sigma^2\Arrowvert C\Arrowvert_{S_2}^2.
\end{equation} 
As compared to (\ref{eq: variance 1}), its variance does not depend on the squared dimension $M^2$ but grows like $rM$, which can be substantially smaller. Moreover, 
$$
\E\Big(\Arrowvert \pi_rX\Arrowvert_{S_2}^2-\sigma^2rM-\Arrowvert C\Arrowvert_{S_2}^2-2\sigma\tr\big(\mathrm{E}'C\big)\Big)^2\ =\ 2\sigma^4rM,
$$
that is, $\sigma^{-1}\big(\Arrowvert \pi_rX\Arrowvert_{S_2}^2-\sigma^2rM-\Arrowvert C\Arrowvert_{S_2}^2\big)$
is approximately centered normal with variance $4\Arrowvert C\Arrowvert_{S_2}^2$ if  $\sigma^2rM=o(1)$ in an asymptotic framework, and $4\Arrowvert C\Arrowvert_{S_2}^2$ is the asymptotic efficiency lower bound (\cite{lama00}). The statistics $\Arrowvert\pi_rX\Arrowvert_{S_2}^2-\sigma^2rM$, however, cannot be used for estimating $\Arrowvert C\Arrowvert_{S_2}^2$ since $\pi_r=\pi_r(C)$ depends on $C$ itself and is unknown a priori. Unfortunately,  Theorem \ref{thm: r=1} and Theorem \ref{thm: lower bounds}  demonstrate  that the same result cannot be shown with
$\Arrowvert \hat{\pi}_rX\Arrowvert_{S_2}^2-\sigma^2rM$ in place of $\Arrowvert \pi_rX\Arrowvert_{S_2}^2-\sigma^2rM$, because the bias $\E\Arrowvert\hat{\pi}_r X\Arrowvert_{S_2}^2-\E\Arrowvert\pi_r X\Arrowvert_{S_2}^2$ is of the order not smaller than $\sigma^2r(M-r)$, i.e. $$\sigma^{-1}\big(\E\Arrowvert\hat{\pi}_r X\Arrowvert_{S_2}^2-\E\Arrowvert\pi_r X\Arrowvert_{S_2}^2\big)$$ 
is not negligible under the same conditions, even not for very favorable matrices $C=\alpha\id_r$. That is, empirical low-rank projections $\Arrowvert\hat{\pi}_rX\Arrowvert_{S_2}^2-\sigma^2rM$  cannot be successively used for efficient estimation of $\Arrowvert C\Arrowvert_{S_2}^2$, even if the $\mathrm{rank}(C)\ll M$ is explicitly known beforehand.
Note that in contrast, empirical low-rank approximations have been proved useful when estimating a low-rank matrix  $C$ under Hilbert-Schmidt norm loss, see \cite{bushwe11}, \cite{kolt12}, \cite{newa11}, and \cite{rots11}. 

The problem of quadratic functional estimation in the matrix context appears, for instance,  in the recent area of low-rank matrix recovery, when one is interested in recovering the linear entropy $1-\Arrowvert C\Arrowvert_{S_2}^2$ of a quantum density matrix $C$ as an approximation of von Neumann entropy based on noisy observations. We refer the reader to \cite{argigu04}  for a detailed description of applications in quantum state tomography and the recent article of \cite{kolt12} for low-rank matrix recovery of quantum density matrices. Note however that our results do not take into account that a quantum density matrix $C$ is self-adjoint, and the Wigner ensemble may behave differently as already outlined in the introduction.  In view of model selection issues, an estimate of the bias is even required over the whole scale $r\in\{1,...,M\}$ since the rank is typically unknown a priori and low at most approximately, for which reason exact asymptotic results for uniformly bounded rank perturbations are of limited value for this application.

%\subsection{Towards high-dimensional principal component analysis} 

\section{Proof of Theorem \ref{thm: general}}\label{sec: theorem main}
This section is devoted to the proof of our main result. Subsection \ref{subsec: slices} deals with the description of the sets $\GG_{M,r}(\delta):=\{\tilde{\pi}_r\in\SS_{M,r}: \Arrowvert \pi_rC\Arrowvert_{S_2}^2 -\Arrowvert\tilde{\pi}_r C\Arrowvert_{S_2}^2\leq\delta\}$, which characterize the slices $\AA_{C,k}= \GG_{M,r}\big(2^{-k}\Arrowvert\pi_r C\Arrowvert_{S_2}^2\big)\setminus \GG_{M,r}\big(2^{-{(k+1)}}\Arrowvert\pi_r C\Arrowvert_{S_2}^2\big)$. It is shown that  these sets can be approximated by sets of  simple geometric structure. In case of some substantial spectral gap in $r$, this approximation is very tight. Sharp bounds on their metric entropy are derived in Subsection \ref{sec: rand}. The final arguments  are given in Subsection \ref{subsec: general}.

\vspace{-4mm}
\subsection{Characterizing the sets $\GG_{M,r}(\delta)$}\label{subsec: slices}
The first goal for a sophisticated analysis is to characterize the sets 
$$
\GG_{M,r}(\delta,C)\ :=\ \Big\{\tilde{\pi}_r\in\SS_{M,r}:\, \Arrowvert \pi_rC\Arrowvert_{S_2}^2 -\Arrowvert\tilde{\pi}_r C\Arrowvert_{S_2}^2\leq\delta\Big\},\ \ \  \delta>0.
$$
Note that $$
\AA_{C,k}\ =\ \GG_{M,r}\Big(2^{-k}\Arrowvert\pi_r C\Arrowvert_{S_2}^2\Big)\setminus \GG_{M,r}\Big(2^{-(k+1)}\Arrowvert\pi_r C\Arrowvert_{S_2}^2\Big).
$$
This is a very delicate part and quite involved in general, but we find below a tight characterization in 
case of approximate rank-$r$-matrices, i.e. those matrices for which $
\Arrowvert(\id-\pi_r)C\Arrowvert_{S_2}$ is small.
For reasons of clarity, we restrict attention to matrices of rank larger or equal to $r$.  %Note that a similar result can be obtained for matrices which can be well approximated by some low-rank projection.

\begin{proposition}\label{prop: 2}
Let $C\in\R^{M\times M}$, $\mathrm{rank}(C)\geq r$, with singular values $\lambda_i,\, i=1,...,M$, ordered in decreasing magnitude. Denote 
$
\Delta_r^*:=\sum_{i=r+1}^{2r}\lambda_i^2$ and $\gamma_r^{*}:=(\lambda_r^2-\lambda_{r+1}^2)^{-1/2}$ if $\lambda_{r}>\lambda_{r+1}$, and $\gamma^*=\infty$ else. 
Then for any $\pi_r\in\mathrm{Arg}\max_{\tilde{\pi}_r\in\SS_{M,r}}\Arrowvert\tilde{\pi}_rC\Arrowvert_{S_2}^2$,
\begin{align}
\GG_{M,r}(\delta,C) \subset &\bigg\{\tilde{\pi}_r\in\SS_{M,r}: \Arrowvert \pi_r -\tilde{\pi}_r\Arrowvert_{S_2}\leq \min\bigg(\lambda_r^{-1}\sqrt{2(\delta +\Delta_r^*)},\, \gamma_r^*\sqrt{2\delta}\bigg)\bigg\}  \label{eq: upper}
\intertext{and}
&  \Big\{\tilde{\pi}_r\in\SS_{M,r}: \Arrowvert \pi_r - \tilde{\pi}_r\Arrowvert_{S_2}\leq \lambda_1^{-1}\sqrt{2\delta}\Big\} \subset \GG_{M,r}(\delta,C).\label{eq: lower}
\end{align}
\end{proposition}

\paragraph{\sc Proof}Let $U\Lambda V'$ denote some singular value decomposition of $C$, with $$\Lambda_r:=\mathrm{diag}\big(\lambda_1,\lambda_2,...,\lambda_r,0,...,0\big)\ \ \text{and}\ \  \Lambda_{M-r}:=\Lambda-\Lambda_r.
$$ 
Recall that $\pi_r:=\sum_{i=1}^rU_iU_i'$ is some maximizer of $\Arrowvert\tilde{\pi}_rC\Arrowvert_{S_2}^2$ over $\tilde{\pi}_r\in\SS_{M,r}$, where the $U_i$'s denote the column vectors of $U$. Note at this point that due to multiplicities in the singular value spectrum, the orthonormal vectors $U_1$, ..., $U_r$ are not unique in general. As concerns the proof of (\ref{eq: upper}), we check first that
\begin{equation}\label{eq: first}
\Arrowvert \pi_r C\Arrowvert_{S_2}^2\ -\ \Arrowvert\tilde{\pi}_r C\Arrowvert_{S_2}^2\  \geq\ \lambda_r^2\frac{1}{2}\Arrowvert  \pi_r-\tilde{\pi}_r\Arrowvert_{S_2}^2-\Xi_r
\end{equation}
for all $\tilde{\pi}_r\in\SS_{M,r}$, where $\Xi_r:=\Arrowvert \tilde{\pi}_rU\Lambda_{M-r}U'\Arrowvert_{S_2}^2=\Arrowvert (\id-\pi_r)\tilde{\pi}_r(\id-\pi_r)U\Lambda_{M-r}U'\Arrowvert_{S_2}^2$.
 Using  the symmetry of the projection matrices $\pi_r$ and $\tilde{\pi}_r$ and the invariance of the trace operator under cyclic permutation, we obtain the identity
\begin{align}
 \Arrowvert \pi_rC&\Arrowvert_{S_2}^2-\Arrowvert \tilde{\pi}_rC\Arrowvert_{S_2}^2\\
&=\ \tr\Big(C'\big(\pi_r'\pi_r-\tilde{\pi}_r'\tilde{\pi}_r\big)C\Big)\nonumber\\
%&=\ \tr\Big(C'\big(\pi_r-\tilde{\pi}_r\big)C\Big)\\
&=\ \tr\Big(CC'\big(\pi_r-\tilde{\pi}_r\big)\Big)\nonumber\\
&=\ \tr\Big(\Lambda_r^2U'\big(\pi_r-\pi_r\tilde{\pi}_r\pi_r\big)U\Big) -  \tr\Big(\Lambda_{M-r}^2U'(\id-\pi_r)\tilde{\pi}_r(\id-\pi_r)U\Big)\nonumber\\
&=\ \sum_{i=1}^r\lambda_i^2\big[U'\big(\pi_r(\id-\tilde{\pi}_r)\pi_r\big)U\big]_{ii}\label{eq: sum1}\\
&\ \ \ \ \ \ \ \ \ \ \ \ \ \ \ \ \  - \sum_{i=r+1}^M\lambda_i^2\big[U'(\id-\pi_r)\tilde{\pi}_r(\id-\pi_r)U\big]_{ii}\label{eq: sum2}
\end{align}
Note that the sum in (\ref{eq: sum2}) equals $\Xi_r$. Since $\pi_r-\pi_r\tilde{\pi}_r\pi_r=\pi_r(\id -\tilde{\pi}_r)\pi_r=\pi_r'(\id-\tilde{\pi}_r)\pi_r$ is positive semidefinite, all summands in the first term of (\ref{eq: sum1}) are non-negative. Consequently, 
\begin{align}
(\ref{eq: sum1})\ &\geq\  \lambda_r^2\sum_{i=1}^r\big[U'\big(\pi_r-\pi_r\tilde{\pi}_r\pi_r\big)U\big]_{ii} \nonumber\\ 
&=\ \lambda_r^2\sum_{i=1}^M\big[U'\big(\pi_r-\pi_r\tilde{\pi}_r\pi_r\big)U\big]_{ii}\label{eq: gleichung}\\
&=\ \lambda_r^2\tr \big(\pi_r-\pi_r\tilde{\pi}_r\pi_r\big),\nonumber
\end{align}
where (\ref{eq: gleichung}) follows from $\big[U'\big(\pi_r-\pi_r\tilde{\pi}_r\pi_r\big)U\big]_{ii}=0$ for $i>r$.  
By positive semidefiniteness again and symmetry, the eigenvalue decomposition of $\pi_r-\pi_r\tilde{\pi}_r\pi_r$ and the invariance of the trace operator under basis transformation yield  $\tr \big(\pi_r-\pi_r\tilde{\pi}_r\pi_r\big)=\Arrowvert \pi_r-\pi_r\tilde{\pi}_r\pi_r\Arrowvert_{S_1}$, hence
$$
 \Arrowvert \pi_r C\Arrowvert_{S_2}^2\ -\ \Arrowvert\tilde{\pi}_r C\Arrowvert_{S_2}^2\ \geq\ \lambda_r^2\Arrowvert \pi_r-\pi_r\tilde{\pi}_r\pi_r\Arrowvert_{S_1} - \Xi_r.
$$
Now, since $\pi_r=\pi_r\tilde{\pi}_r+\pi_r(\id-\tilde{\pi}_r)$ is an orthogonal decomposition, $$\Arrowvert \pi_r\tilde{\pi}_r\Arrowvert_{S_2}^2+\Arrowvert\pi_r(\id-\tilde{\pi}_r)\Arrowvert_{S_2}^2 = \Arrowvert\pi_r\Arrowvert_{S_2}^2=r=\Arrowvert\tilde{\pi}_r\Arrowvert_{S_2}^2=\Arrowvert \pi_r\tilde{\pi}_r\Arrowvert_{S_2}^2+\Arrowvert (\id-\pi_r)\tilde{\pi}_r\Arrowvert_{S_2}^2,
$$ 
implying that  $\Arrowvert \pi_r(\id -\tilde{\pi}_r)\Arrowvert_{S_2}^2=\Arrowvert (\id -\pi_r)\tilde{\pi}_r\Arrowvert_{S_2}^2$. Consequently,
\begin{align}
\Arrowvert \pi_r-\pi_r\tilde{\pi}_r\pi_r\Arrowvert_{S_1}=\tr \big(\pi_r-\pi_r\tilde{\pi}_r\pi_r\big)\ &=\ \tr \big(\pi_r'\pi_r-\pi_r'\tilde{\pi}_r'\tilde{\pi}_r\pi_r\big)\label{eq: chain}\\
&=\ \Arrowvert\pi_r\Arrowvert_{S_2}^2-\Arrowvert\pi_r\tilde{\pi}_r\Arrowvert_{S_2}^2\nonumber\\
&=\ \Arrowvert\pi_r(\id-\tilde{\pi}_r)\Arrowvert_{S_2}^2\ =\ \frac{1}{2}\Arrowvert \pi_r-\tilde{\pi}_r\Arrowvert_{S_2}^2,\nonumber
\end{align}
where the last equality follows by the orthogonality of the decomposition $\tilde{\pi}_r-{\pi}_r=(\id-\pi_r)\tilde{\pi}_r-\pi_r(\id-\tilde{\pi}_r)$. This implies (\ref{eq: first}). In order to deduce the  bound in (\ref{eq: upper}), note first that because of the positive semidefiniteness of $(\id-\pi_r)\tilde{\pi}_r(\id-\pi_r)$, also all summands in the second sum (\ref{eq: sum2}) are non-negative, whence
\begin{align}
 \sum_{i=r+1}^M\lambda_i^2\big[U(\id-\pi_r) \tilde{\pi}_r(\id-\pi_r)U'\big]_{ii}\ \leq\ \lambda_{r+1}^2\sum_{i=r+1}^M\big[U(\id-\pi_r) \tilde{\pi}_r(\id-\pi_r)U'\big]_{ii}.\nonumber%\\
%&\leq\ \lambda_{r+1}^2\frac{1}{2}\Arrowvert \pi_r-\tilde{\pi}_r\Arrowvert_{S_2}^2\label{eq: first1}
\end{align}
With the same arguments as provided  above for $\pi_r-\pi_r\tilde{\pi}_r\pi_r=\pi_r(\id-\tilde{\pi}_r)\pi_r$, we deduce
\begin{align}
\Xi_r\ \leq\ 
 \lambda_{r+1}^2\frac{1}{2}\Arrowvert (\id-\pi_r)-(\id-\tilde{\pi}_r)\Arrowvert_{S_2}^2\ =\ \lambda_{r+1}^2\frac{1}{2}\Arrowvert \pi_r-\tilde{\pi}_r\Arrowvert_{S_2}^2.\label{eq: first1}
\end{align} Moreover, because $U'\tilde{\pi}_rU\in\SS_{M,r}$ for every $U\in O(M)$,
\begin{equation}\Xi_r\  =\  \Arrowvert \tilde{\pi}_rU\Lambda_{M-r}U'\Arrowvert_{S_2}^2\ =\ \Arrowvert U'\tilde{\pi}_rU\Lambda_{M-r}\Arrowvert_{S_2}^2\  \leq\ \sum_{i=r+1}^{2r}\lambda_i^2,\label{eq: first2}
\end{equation}
and claim (\ref{eq: upper}) follows from (\ref{eq: first}), together with (\ref{eq: first1}) and (\ref{eq: first2}).  The proof of (\ref{eq: lower}) uses that the expression  (\ref{eq: sum1}) is in turn upper bounded by $\lambda_1^2\Arrowvert \pi_r-\pi_r\tilde{\pi}_r\pi_r\Arrowvert_{S_1}$ while (\ref{eq: sum2}) is less or equal to zero, and concludes finally with the same equality chain (\ref{eq: chain}).\hfill $\square$

We note that in the particular case of rank-$r$-matrices, $\sum_{i=r+1}^{2r}\lambda_i=0$, i.e. the first term in the upper bound on the radius of the $S_2$-ball in (\ref{eq: upper})  coincides up to the ratio $\lambda_1/\lambda_r$ with the lower bound in (\ref{eq: lower}). Equality holds for rank-$r$-matrices with rectangular singular value spectrum $\lambda_1=\lambda_2=\cdots=\lambda_r$. If $C=\alpha\id$ for some $\alpha\not= 0$, then $\Delta_r^*=\alpha^2r$ and the inclusion (\ref{eq: upper}) is trivial: 
$$
\bigg\{\tilde{\pi}_r\in\SS_{M,r}: \Arrowvert \pi_r -\tilde{\pi}_r\Arrowvert_{S_2}\leq \alpha^{-1}\sqrt{2(\delta +\alpha^2r)}\bigg\}\ =\ \SS_{M,r},\ \ \text{for any $\delta>0$}.
$$
This is in accordance with the fact that in this case, also $\GG_{M,r}(\delta)=\SS_{M,r}$ for any $\delta>0$.

\subsection{Metric entropy bounds}\label{sec: rand}
The nice feature of the results from the previous paragraph is that they enable us to determine tight bounds on the metric entropy on these particular subsets of the Grassmannian by the volumetric argument. We provide a slightly refined version. Recall at this point the definition of the covering numbers. For any totally-bounded, pseudometric space $(\mathcal{X},d)$ and any subset $E\subset\XX$, the covering number $N(E,d,\delta)$ is the smallest number of closed $d$-balls in $\XX$ of radius $\delta$ needed to cover $E$.

\begin{lemma}\label{lemma: bounds}
For any $\pi_r\in\SS_{M,r}$ and $\delta>0$, let $B_{S_q}(\pi_r,\delta)$ denote the closed $S_q$-ball with center $\pi_r$ of radius $\delta$. Then there exist universal constants $\hat{c}, C, c>0$ such that for all $0<\Delta\leq r$, $r\leq M-r$ and $0<\delta\leq \Delta$,
\begin{align}
\log N\Big(B_{S_2}(\pi_r,\Delta)\cap\SS_{M,r}, d_{S_2},\delta\Big)
&\leq \min\hspace{-0.5mm}\left(\hat{c}rM\frac{\Delta^2}{\delta^2},r(M-r)\hspace{-0.5mm}\log\hspace{-0.5mm}\bigg(\frac{C\Delta}{\delta}\bigg)\hspace{-0.8mm}\right),\label{eq: S2}\\
\log N\Big(B_{S_2}(\pi_r,\Delta)\cap\SS_{M,r}, d_{S_{\infty}},\delta\Big) &\leq  \min\hspace{-0.5mm}\left(\hat{c}rM\frac{\Delta}{\sqrt{r}\delta}, r(M-r)\hspace{-0.5mm}\log\hspace{-0.5mm}\bigg(\frac{C\Delta}{\delta}\bigg)\hspace{-0.8mm}\right)\label{eq: Sinfty}
\end{align}
as well as

\vspace{-9mm}
\begin{align}
\log N\Big(B_{S_2}(\pi_r,\Delta)\cap \SS_{M,r}, d_{S_2},\delta\Big)\ &\geq\ r(M-r)\log\bigg(\frac{c\Delta}{\delta}\bigg)\ \ \ \text{and}\label{eq: S2-2}\\
\log\hspace{-0.5mm} N\Big(B_{S_2}(\pi_r,\Delta)\hspace{-0.5mm}\cap\hspace{-0.5mm}\SS_{M,r}, d_{S_{\infty}},\delta \Big)\ &\geq\ r(M-r)\log\bigg(\frac{c\Delta}{\sqrt{r}\delta}\bigg)\label{eq: Sinfty-2}
\end{align}
\end{lemma}

\vspace{-6mm}
\paragraph{\sc Remark} (i) $\hat{c}^{1/2}$ equals $c_{\mathrm{Sud}}$ times a uniform (in $M$) bound on the expected spectral norm of $\mathrm{E}/(\sqrt{M}\sigma)$, where $c_{\mathrm{Sud}}$ is the universal constant of Sudakov's minoration which is bounded by $6$ (see \cite{led96}). $C$ is proportional to the ratio of the constants $C'/c'$ which appear in the bounds of the metric entropy of the Grassmann manifold, see the proof below. Estimates on their values are not provided in \cite{sza82}. 

(ii) For arbitrary $1\leq q\leq\infty$, the bound
\begin{equation}\label{eq: S-infty}
r(M-r)\log\bigg(\frac{c\Delta}{\delta}\bigg) \ \leq\ \log N\Big(B_{S_q}(\pi_r,\Delta)\cap\SS_{M,r}, d_{S_q},\delta\Big)\
\leq \ r(M-r)\log\bigg(\frac{C\Delta}{\delta}\bigg)
\end{equation}
can be proved completely analogously, replacing below $S_2$ by $S_q$.
\vspace{-3mm}
\paragraph{\sc Proof}
By the geometric formulation of Sudakov's minoration, the trace duality and the Cauchy-Schwarz inequaliy,
\begin{align*}
\delta\sqrt{\log N\Big(B_{S_2}(\pi_r,\Delta)\cap \SS_{M,r}\, d_{S_2},\delta\Big)}\ &\lesssim\ \E\sup_{T\in B_{S_2}(\pi_r,\Delta)\cap \SS_{M,r}}\tr\big(T(\mathrm{E}/\sigma)\big)\\
&\leq\ \E\sup_{T\in B_{S_2}(\pi_r,\Delta)\cap \SS_{M,r}}\tr\big((T-\pi_r)(\mathrm{E}/\sigma)\big) + \E\,\tr\big(\pi_r(\mathrm{E}/\sigma)\big)\\
&\leq\ \sup_{T\in B_{S_2}(\pi_r,\Delta)\cap \SS_{M,r}}\Arrowvert T-\pi_r\Arrowvert_{S_1} \E\Arrowvert \mathrm{E}/\sigma\Arrowvert_{S_{\infty}}\\
&\lesssim\ {\Delta}\sqrt{rM},
\end{align*}
which provides the first estimate in the minimum of (\ref{eq: S2}). %We say that a set $\FF_{\gamma}\subset\SS_{M,r}$ is a $\gamma$-net for $\SS_{M,r}$ with respect to some metric $d$, if 
%$
%\SS_{M,r}\subset\cup_{\tilde{\pi}_r\in\FF_{\gamma}} B(\tilde{\pi}_r,d,\gamma),
%$
%where $B(\tilde{\pi}_r,d,\gamma)$ denotes the closed $d$-ball with center $\tilde{\pi}_r$ and radius $\gamma$. For $\Delta>\delta$, let $\FF_{\delta}$  be a minimal $\delta$-net with respect to $d_{S_2}$ in $\SS_{M,r}$. Then
%\begin{align}
%N\Big(B_{S_2}(\pi_r,\Delta)\cap \SS_{M,r},d_{S_2},\delta\Big)\ &\leq\ \mathrm{card}\Big\{\tilde{\pi}_r\in\FF_{\delta}:\, \tilde{\pi}_r\in B_{S_2}\big(\pi_r,2\Delta\big)\Big\}\nonumber\\
%&\leq\ N\Big(B_{S_2}(\pi_r,4\Delta)\cap \SS_{M,r},d_{S_2},\delta/2\Big),\nonumber
%\end{align}
%where the second inequality follows from the minimality of $\FF_{\delta}$.
%Since $\Arrowvert \tilde{\pi}_r^{(1)}-\tilde{\pi}_r^{(2)}\Arrowvert_{S_{\infty}}\leq 2$ and consequently $\Arrowvert \tilde{\pi}_r^{(1)}-\tilde{\pi}_r^{(2)}\Arrowvert_{S_2}\leq 2\Arrowvert \tilde{\pi}_r^{(1)}-\tilde{\pi}_r^{(2)}\Arrowvert_{S_1}$ for any two elements $ \tilde{\pi}_r^{(1)},\tilde{\pi}_r^{(2)}\in\SS_{M,r}$, 
In order to prove the second term, note first that $N\big(B_{S_2}(\pi_r,\Delta)\cap \SS_{M,r},d_{S_2},\delta\big)$ does not depend on the specific $\pi_r\in\SS_{M,r}$. Similarly to the covering number, the capacity number $D(E,d,\delta)$ is the largest number of elements of $E$ having distance $d$ strictly larger than $\delta$ to each other. Using the relations between covering and packing (capacity) 
numbers (cf.~Theorem 1.2.1, \cite{dud99}), 

\vspace{-8mm}
\begin{align*}
N\Big(B_{S_2}(\pi_r,\Delta)&\cap \SS_{M,r},d_{S_2},\delta\Big)N\big(S_{M,r}, d_{S_2}, 4\Delta\big)\\ 
&\leq\ D\Big(B_{S_2}(\pi_r,\Delta)\cap \SS_{M,r},d_{S_2},\delta\Big)D\big(S_{M,r}, d_{S_2}, 4\Delta\big).
\end{align*}
Let $\{\pi_r^{(1)},...,\pi_r^{(k_{4\Delta})}\}$ be some maximal subset of $\SS_{M,r}$ with $d_{S_2}(\pi_r^{(l)},\pi_r^{(m)})>4\Delta$ for all $l\not=m\in\{1,...,k_{4\Delta}\}$,  $k_{4\Delta}=D(\SS_{M,r},d_{S_2},4\Delta)$. Then 
$$
\sum_{j=1}^{k_{4\Delta}}D\big(B_{S_2}(\pi_r^{(j)},\Delta)\cap \SS_{M,r},d_{S_2},\delta\big)\ \leq\  D(\SS_{M,r},d_{S_2},\delta/2)\ \leq\ N\big(\SS_{M,r},d_{S_2},\delta/4\big),
$$
that is
$~~~
N\Big(B_{S_2}(\pi_r,\Delta)\cap \SS_{M,r},d_{S_2},\delta\Big)\ \leq\ N\big(\SS_{M,r},d_{S_2},\delta/4\big)\Big/N\big(S_{M,r}, d_{S_2}, 4\Delta\big)$.

Now (\ref{eq: S2}) follows by an immediate application of Proposition 8, \cite{pajo98}, which states
$$
\Big(\frac{c'}{\xi}\Big)^{r(M-r)}\ \leq\ N\big(\SS_{M,r}, d_{S_q},\xi r^{1/q}\big)\ \leq\ \Big(\frac{C'}{\xi}\Big)^{r(M-r)}
$$
for $1\leq q\leq\infty$ and universal constants $c',C'>0$. As concerns bound (\ref{eq: Sinfty}), first observe that
\begin{equation}
 N\Big(B_{S_2}(\pi_r,\Delta)\cap\SS_{M,r}, d_{S_{\infty}},\delta\Big) \leq  N\Big(B_{S_2}(\pi_r,\Delta)\cap\SS_{M,r}, d_{S_2},\delta/\theta\Big) N\Big(B_{S_2}(0,1), d_{S_{\infty}},\theta\Big)\label{eq: factor}
\end{equation}
for any $\theta>0$. Combining (\ref{eq: factor}) with (\ref{eq: S2}) and the bound
$
\log N\big(B_{S_2}(0,1), d_{S_{\infty}},\delta\big) \lesssim M\delta^{-2}
$
(cf.~Pajor 1998, Lemma 4), we obtain
$$
\log N\Big(B_{S_2}(\pi_r,\Delta)\cap\SS_{M,r}, d_{S_{\infty}},\delta\Big)\ \lesssim\ rM\frac{\Delta^2\theta^2}{\delta^2}+M\theta^{-2}.
$$
Choosing $\theta^2=\delta\Delta^{-1}r^{-1/2}$ gives the first term in the minimum on the RHS of (\ref{eq: Sinfty}). The proof of the second bound in (\ref{eq: Sinfty}) follows from (\ref{eq: S2}) since $d_{S_{\infty}}\leq d_{S_{2}}$. As concerns the reverse inequality (\ref{eq: S2-2}), it is sufficient to note that
$$
D\big(\SS_{M,r},d_{S_2},\Delta\big)N\Big(B_{S_2}(\pi_r,\Delta)\cap \SS_{M,r},d_{S_2},\delta\Big)\ \geq\ N\big(\SS_{M,r},d_{S_2},\delta\big)%\Big(B_{S_2}(\pi_r,2\Delta)\cap \SS_{M,r},d_{S_2},\delta\Big)% \mathrm{card}\Big\{\tilde{\pi}_r\in\FF_{\delta}:\, \tilde{\pi}_r\in B_{S_2}\big(\pi_r,\Delta\big)\Big\},
$$
which after applying the inequalities of Theorem 1.2.1 in \cite{dud99} again is lower bounded by
\begin{equation}\nonumber
%\inf_{\pi_r\in\SS_{M,r}}\mathrm{card}\Big\{\tilde{\pi}_r\in\FF_{\delta}:\, \tilde{\pi}_r\in B_{S_2}\big(\pi_r,\Delta\big)\Big\}
N\Big(B_{S_2}(\pi_r,\Delta)\cap \SS_{M,r},d_{S_2},\delta\Big)\ \geq\ \frac{N\big(\SS_{M,r},d_{S_2},\delta\big)}{N\big(S_{M,r}, d_{S_2}, \Delta/2\big)},
\end{equation}
and the result follows as above by an application of Proposition 8, \cite{pajo98}. (\ref{eq: Sinfty-2}) follows analogously from $d_{S_{\infty}}\geq (2r)^{-1/2}d_{S_{2}}$ on $\SS_{M,r}$.  \hfill$\square$

\subsection{Slicing the Grassmann manifold}\label{subsec: general}
As has been seen in  Section \ref{sec: lower bounds}, the bound involves some term of the order $\sigma^2r(M-r)$ at least. Since Proposition \ref{prop: rough bound} yields in case $C=0$ the bound $\delta_{0,M,\sigma^2,r}\lesssim \sigma^2rM$, we decompose the supremum
\begin{align*}
\sup_{\tilde{\pi}_r\in\SS_{M,r}}&\Big(\big\Arrowvert \hat{\pi}_r (C+\mathrm{E})\big\Arrowvert_{S_2}^2  - \big\Arrowvert \pi_r ( C+\mathrm{E})\big\Arrowvert_{S_2}^2\Big)\\ 
&=\ \sup_{\tilde{\pi}_r\in\SS_{M,r}}\bigg\{\big\Arrowvert\tilde{\pi}_r\mathrm{E}\big\Arrowvert_{S_2}^2-\big\Arrowvert{\pi}_r\mathrm{E}\big\Arrowvert_{S_2}^2 +2\, \tr\Big(\mathrm{E}'(\tilde{\pi}_r-\pi_r)C\Big)\\
&\ \ \ \ \ \ \ \ \ \ \ \ \ \ \ \ \ \ \ \ \ \ \  \ \ \ \ \ \ \ \ \ \ \ \ \ \ \ \ \ \ \ \ \ \ \ \ \ \ \ \ \ {-\Big(\big\Arrowvert{\pi}_rC\big\Arrowvert_{S_2}^2-\big\Arrowvert \tilde{\pi}_rC\big\Arrowvert_{S_2}^2\Big)}\bigg\}\\
&\leq\ \sup_{\tilde{\pi}_r\in\SS_{M,r}}\Big\{\big\Arrowvert\tilde{\pi}_r\mathrm{E}\big\Arrowvert_{S_2}^2-\big\Arrowvert{\pi}_r\mathrm{E}\big\Arrowvert_{S_2}^2\Big\}\\
&\ \ \ \ \ \ \ \ \ \ \ \  +\ \sup_{\tilde{\pi}_r\in\SS_{M,r}}\bigg\{2 \tr\Big(\mathrm{E}'(\tilde{\pi}_r-\pi_r)C\Big) {-\Big(\big\Arrowvert{\pi}_rC\big\Arrowvert_{S_2}^2-\big\Arrowvert \tilde{\pi}_rC\big\Arrowvert_{S_2}^2\Big)}\bigg\}
\end{align*}
and treat these two suprema separately. Since
$$
\E\bigg( \sup_{\tilde{\pi}_r\in\SS_{M,r}}\Big\{\big\Arrowvert\tilde{\pi}_r\mathrm{E}\big\Arrowvert_{S_2}^2-\big\Arrowvert{\pi}_r\mathrm{E}\big\Arrowvert_{S_2}^2\Big\}\bigg)\ \lesssim\ \sigma^2rM
$$
by Proposition \ref{prop: rough bound} applied to the situation $C=0$, 
it is sufficient to prove that the expectation of
\begin{align*}
W\ :=&\ \sup_{\tilde{\pi}_r\in\SS_{M,r}}\bigg\{2 \tr\Big(\mathrm{E}'(\tilde{\pi}_r-\pi_r)C\Big) {\ -\ \Big(\big\Arrowvert{\pi}_rC\big\Arrowvert_{S_2}^2-\big\Arrowvert \tilde{\pi}_rC\big\Arrowvert_{S_2}^2\Big)}\bigg\}
\end{align*}
satisfies the bound of Theorem \ref{thm: general}. For any $C\in\R^{M\times M}$, $\mathrm{rank}(C)\geq r$, let $C=U\Lambda V'$ denote some singular value decomposition, with diagonal matrix  $\Lambda=\mathrm{diag}(\lambda_1,\lambda_2,...,\lambda_M)$. The singular values are assumed to be ordered in decreasing magnitude. Let 
$$\SS_{M,r}(\delta)\ :=\ \Big\{\tilde{\pi}_r\in\SS_{M,r}: \Arrowvert \tilde{\pi}_r-\pi_r\Arrowvert_{S_2}\leq \sqrt{2\delta}\Big\}.
$$
In view of the inclusion provided by Proposition \ref{prop: 2}, we deviate slightly from the description at the beginning of Section \ref{subsec: iteration} and conduct the proof of Theorem \ref{thm: general} along two different decompositions of $\SS_{M,r}$. We shall decompose $\SS_{M,r}$ into slices $$\CC_{C,i}(\tilde{\Delta}_k)\ =\ \BB_{C,i}(\tilde{\Delta}_k)\setminus\BB_{C,i}(\tilde{\Delta}_{k-1})
$$ 
along a geometric grid $\tilde{\Delta}_k= 2^{-k+2}r$ up to $k< k_{0,i}$ with $k_{0,i}$ specified below. First, we take
\begin{equation}
\BB_{C,1}(\tilde{\Delta}_k)\ := \ \Big\{\tilde{\pi}_r\in\SS_{M,r}:\, \Arrowvert\pi_rC\Arrowvert_{S_2}^2-\Arrowvert\tilde{\pi}_rC\Arrowvert_{S_2}^2\geq \tilde{\Delta}_k\lambda_r^2-\sum_{i=r+1}^{2r}\lambda_i^2\Big\}.
\end{equation}
In a second step, whenever $\lambda_{r}>\lambda_{r+1}$, we choose
\begin{equation}
\BB_{C,2}(\tilde{\Delta}_k)\ := \ \Big\{\tilde{\pi}_r\in\SS_{M,r}:\, \Arrowvert\pi_rC\Arrowvert_{S_2}^2-\Arrowvert\tilde{\pi}_rC\Arrowvert_{S_2}^2\geq (\lambda_r^2-\lambda_{r+1}^2) \tilde{\Delta}_k\Big\}.
\end{equation}
By Proposition \ref{prop: 2}, $\CC_{C,i}(\tilde{\Delta}_k)\subset\SS_{M,r}(2\tilde{\Delta}_k)$.
Recall that by construction,
$
\Arrowvert \pi_rC\Arrowvert_{S_2}^2 - \Arrowvert \tilde{\pi}_rC\Arrowvert_{S_2}^2  \geq   \tilde{\Delta}_k\lambda_r^2 - \sum_{i=r+1}^{2r}\lambda_i^2$ if $\tilde{\pi}_r\in\CC_{C,1}(\tilde{\Delta}_k).
$ 
Whenever $\lambda_{r+1}>0$, define
\begin{equation}\label{eq: k0}
k_{0,1}\ :=\ \arg\max_{k\in\N_0}\left\{\tilde{\Delta}_k\lambda_r^2 - \sum_{i=r+1}^{2r}\lambda_i^2\geq\frac{1}{2}\tilde{\Delta}_k\lambda_r^2\right\},\ \  \text{and set}\  k_{0,1}:=\infty\ \text{if}\ \lambda_{r+1}=0.
\end{equation}
Define $k_{0,2}:=\infty$. Denoting
\begin{align*}
W_{k,i}\ :=&\ \sup_{\tilde{\pi}_r\in\CC_{C,i}(\tilde{\Delta}_k)}\bigg\{2\, \tr\Big(\mathrm{E}'(\tilde{\pi}_r-\pi_r)C\Big) {\ -\ \Big(\big\Arrowvert{\pi}_rC\big\Arrowvert_{S_2}^2-\big\Arrowvert \tilde{\pi}_rC\big\Arrowvert_{S_2}^2\Big)}\bigg\},\\
W_{k,i}^0\ :=&\ \sup_{\tilde{\pi}_r\in\CC_{C,i}(\tilde{\Delta}_{k})}2\, \tr\Big(\mathrm{E}'(\tilde{\pi}_r-\pi_r)C\Big)
\end{align*} 
and $$\tilde{W}_{k_{0,1}}:=\sup_{\tilde{\pi}_r\in\SS_{M,r}(\tilde{\Delta}_{k_{0,1}})} 2\, \tr\Big(\mathrm{E}'(\tilde{\pi}_r-\pi_r)C\Big)$$
with $\tilde{\Delta}_{\infty}:=0$, we obtain the expansions
\begin{align}
\E\,W\
&\leq\  \begin{cases}\sum_{k < k_{0,1}}\E\underset{~}{\Big(}0\vee \big(W_{k,1}^0-\Omega_{k,1}\big)\Big)\  +\ \E\,\tilde{W}_{k_{0,1}}& \\
\sum_{k \in \N}\E\overset{~}{\Big(}0\vee \big(W_{k,2}^0-\Omega_{k,2}\big)\Big) & 
\end{cases}\label{eq: expansion-m}
\end{align}
where $\Omega_{k,1}=\tilde{\Delta}_k\lambda_r^2-\sum_{i=r+1}^{2r}\lambda_i^2$ and $\Omega_{k,2}=(\lambda_r^2-\lambda_{r+1}^2)\tilde{\Delta}_k$. Note that $\E\,\tilde{W}_{\infty}=0$, that is, $\E\,\tilde{W}_{k_{0,1}}=0$ if $\lambda_{r+1}=0$.

\paragraph{\sc Proof of Theorem \ref{thm: general}} Note that each $W_{k,i}^0$ is the supremum over some Gaussian process, and
$$
\sup_{\tilde{\pi}_r\in\CC_{C,i}(\tilde{\Delta}_k)}\mathrm{Var}\bigg(\tr\Big(\mathrm{E}'(\tilde{\pi}_r-\pi_r)C\Big)\bigg)\ \leq \sup_{\tilde{\pi}_r\in\CC_{C,i}(\tilde{\Delta}_k)}\sigma^2\lambda_1^2\Arrowvert \tilde{\pi}_r-\pi_r\Arrowvert_{S_2}^2\ \leq\ 2\sigma^2\lambda_1^2\tilde{\Delta}_k\ =:\ \sigma_k^2,
$$
where the last equality follows from the inclusion $\CC_{C,i}(\tilde{\Delta}_k)\subset\SS_{M,r}(2\tilde{\Delta}_k)$.
By Lemma \ref{lemma: bounds} and $S_2$-chaining over $\SS_{M,r}(2\tilde{\Delta}_k)$, 
\begin{align*}
\E\,\big(0\vee W_{k,i}^0\big)\ &\lesssim\ \int_0^{1}\sigma\lambda_1\tilde{\Delta}_k^{1/2} \Big(2\log\big(N\big(\SS_{M,r}(\tilde{\Delta}_k),d_{S_2}, 2\tilde{\Delta}_k^{1/2}\delta\big)\big)\Big)^{1/2}\,\d \delta\\ 
&\lesssim\ \sigma\lambda_1\tilde{\Delta}_k^{1/2}\sqrt{r(M-r)}. 
\end{align*}
I.e., there exists some constant $c>0$, independent of $M,r,\sigma^2$ and $C$, such that 
\begin{equation}\label{eq: expect}
\E\,\big(0\vee W_{k,i}^0\big)\ \leq\ c\sigma\lambda_1\tilde{\Delta}_k^{1/2}\sqrt{r(M-r)},
\end{equation} and by the \cite{bor75} - \cite{suci74} inequality,
\begin{equation}\label{eq: conc bound}
\P\bigg(W_{k,i}^0 \, \geq\, c\sigma\lambda_1\tilde{\Delta}_k^{1/2}\sqrt{r(M-r)}\, +\, \sqrt{2}\sigma\lambda_1\tilde{\Delta}_k^{1/2}\sqrt{2\eta}\bigg)\ \leq\ \exp(-\eta)\ \ \text{for any $\eta>0$}.
\end{equation}
{\sc (The case $i=1$)} With the help of (\ref{eq: conc bound}) we evaluate the first term 
$$
\sum_{k < k_{0,1}}\E\Big(0\vee \big(W_{k,1}^0-\Omega_{k,1}\big)\Big)\ \leq\ \sum_{k < k_{0,1}}\E\Big(0\vee \big(W_{k,1}^0-\frac{1}{2}\tilde{\Delta}_k\lambda_r^2\big)\Big)$$ 
in expansion (\ref{eq: expansion-m}). For the ease of notation, define $\Delta_k:= c\sigma\lambda_1\tilde{\Delta}_k^{1/2}\sqrt{r(M-r)}$. Since $\Delta_k$ depends only linearly on $\tilde{\Delta}_k^{1/2}$, while $$\Arrowvert\pi_rC\Arrowvert_{S_2}^2-\Arrowvert\tilde{\pi}_rC\Arrowvert_{S_2}^2 \ \geq\ \frac{1}{2}\tilde{\Delta}_k\lambda_r^2\ \ \ \text{for all~}\tilde{\pi}_r\in\CC_C(\tilde{\Delta}_k)\ \text{and}\ k\leq k_{0,1},
$$ 
define the additional auxiliary integer
$$
k_1^*\ :=\ \arg\max_{k \leq k_{0,1}}\left\{\frac{1}{2}\tilde{\Delta}_k\lambda_r^2-\Delta_k\geq \frac{1}{4}\tilde{\Delta}_k\lambda_r^2\right\},
$$
and set $k_1^*=0$ if the relation never holds. If $k_1^*=0$, then 
$$
\sqrt{r}\lambda_r^2\ \leq\ 4c\sigma\lambda_1\sqrt{r(M-r)},$$
and the bound
$$
\E\sup_{\tilde{\pi}_r\in\SS_{M,r}}\Big(\Arrowvert \tilde{\pi}_r(C+\mathrm{E})\Arrowvert_{S_2}^2 -\ \Arrowvert \pi_r(C+\mathrm{E})\Arrowvert_{S_2}^2\Big)\
 \lesssim\ \sigma^2rM\bigg(1+\frac{\lambda_1}{\sigma\sqrt{M}}\bigg)
$$
follows immediately by Corollary \ref{cor: upper} (i). Thus, we assume $k_1^*\geq 1$. We first treat the case $k\leq k_1^*$.
 By the representation formula for the expectation of non-negative random variables and the definition of $k_1^*$,
\begin{align*}
\E\Big(0\vee \big(W_{k,1}^0-\frac{1}{2}\tilde{\Delta}_k\lambda_r^2\big)\Big)\ &=\ \int_0^{\infty}\P\Big(W_{k,1}^0-\frac{1}{2}\tilde{\Delta}_k\lambda_r^2>u\Big)\d u\\
&\leq\ \int_0^{\infty}\P\bigg(Z_k^0-\Delta_k > u + \frac{1}{4}\tilde{\Delta}_k\lambda_r^2\bigg)\d u
\end{align*} 
for all $k\leq k_1^*$.
Next, by (\ref{eq: conc bound}) and with $
A_k = (2\sigma_k)^{-1}\frac{1}{4}\tilde{\Delta}_k\lambda_r^2,
$
\begin{align}
\int_0^{\infty}\P\bigg(W_{k,1}^0>u + \frac{1}{4}\tilde{\Delta}_k\lambda_r^2\bigg)\d u\ \nonumber
&\leq\ \int_0^{\infty}\exp\left(-\frac{\big(u+ \frac{1}{4}\tilde{\Delta}_k\lambda_r^2\big)^2}{4\sigma_k^2}\right)\d u\nonumber\\
%&=\ 2\sqrt{\tilde{\Delta}_k\big(\sigma^4rM+\frac{r\sigma^2}{2}\Arrowvert C\Arrowvert_{S_{\infty}}^2\big)}\int_{A_k}\exp\Big(\frac{u^2}{2}\Big)\d u\\
&\lesssim\ \frac{\sigma_k}{2+A_k}\exp\big(-A_k^2/2\big)\label{eq: integral-m} \\
&\leq\ \frac{\sigma_k}{(2+A_k)(1+A_k^2/2)},\label{eq: integral-m2}
\end{align}
where we used $\P(N\geq x)\leq (2+x)^{-1}\exp(-x^2/2)$ for $N\sim\NN(0,1)$ in (\ref{eq: integral-m}) and  the inequality $\exp(-x)\leq 1/(1+x)$ $\forall\,x>0$ in (\ref{eq: integral-m2}). Thus
\begin{align}
\sum_{k\leq k_1^*}\E\Big(0\vee \big(W_{k,1}^0-\Delta_k-\frac{1}{2}\tilde{\Delta}_k\lambda_r^2\big)\Big)\ \lesssim\ \sum_{k\leq k_1^*}\bigg(\frac{\sigma_k}{1+A_k^2/2}\bigg).\label{eq: k_0}
\end{align}
We evaluate (\ref{eq: k_0}). Recall $\sigma_k=\sqrt{2}\sigma\lambda_1\tilde{\Delta}_k^{1/2}$. Then,  bounding the sum by an intgral and by change of variables,
\begin{align*}
\sum_{k\leq k_1^*}\frac{\sigma_k}{1+A_k^2/2}\ &\leq\ \sum_{k\in\N}\sqrt{2}\sigma\lambda_1\tilde{\Delta}_k^{1/2}\bigg(1+\frac{\tilde{\Delta}_k^2\lambda_r^4}{ 128\sigma^2\lambda_1^2\tilde{\Delta}_k}\bigg)^{-1}\\
&\lesssim\ \sigma{\lambda}_1\int_0^1\bigg(1+ x^2\frac{{\lambda}_r^4}{{\sigma^2\lambda}_1^2}\bigg)^{-1}\d x\\
&\lesssim\  \sigma^2\,\frac{{\lambda}_1^2}{{\lambda}_r^2}\int_0^{\infty}(1+x^2)^{-1}\d x.%\ \lesssim \ \sigma^2\sqrt{rM}\,\frac{1+\tilde{\lambda}_1}{\tilde{\lambda}_r}.
\end{align*}
In order to estimate the expression $\sum_{k\geq k_1^*+1}^{k_{0,1}-1}\E(0\vee (W_{k,1}^0-(1/2)\tilde{\Delta}_k\lambda_r^2))$ we need to determine a lower bound on $k_1^*$ in dependence of $\lambda_1,\lambda_r$.  By its definition, $k > k_1^*$ implies
\begin{equation}\label{eq: k*}
c\sigma\lambda_1\tilde{\Delta}_k^{1/2}\sqrt{r(M-r)}>\tilde{\Delta}_k\lambda_r^2/4
\end{equation}
as long as $k_1^*<k_{0,1}$. Recalling by (\ref{eq: expect}) that $\E(0\vee W_{k,1}^0)\leq \Delta_k$, and using (\ref{eq: k*}) and the representation for the tail of the geometric series,  
\begin{align*}
\sum_{k>k_1^*}^{k_{0,1}-1}\E \big(0\vee W_{k,1}^0\big)\ &\lesssim\ \sum_{k>k_1^*}\tilde{\Delta}_k^{1/2}\sigma\lambda_1\sqrt{r(M-r)}\ \lesssim \ \sigma^2rM\frac{\lambda_1^2}{{\lambda}_r^2}
\end{align*}
As concerns the second term $\E\tilde{W}_{k_{0,1}}$ in expansion (\ref{eq: expansion-m}), we obtain by definition (\ref{eq: k0}) of $k_{0,1}$, if $k_{0,1}<\infty$,
\begin{align*}
\E\,\tilde{W}_{k_{0,1}}\ \lesssim \ \tilde{\Delta}_{k_{0,1}}\sigma\lambda_1\sqrt{r(M-r)}\ \lesssim\ \bigg(\frac{\frac{1}{r}\sum_{i=r+1}^{2r}\lambda_i^2}{\lambda_r^2}\bigg)^{1/2}\cdot\sigma\lambda_1r\sqrt{(M-r)}.
\end{align*}
If $k_{0,1}=\infty$, then $\mathrm{rank}(C)=r$ and the rank-$r$-projection $\pi_r$ is unique, i.e.~$\tilde{W}_{k_{0,1}}=0$. Collecting things together, this proves the bound
\begin{equation}\label{eq: general1}
\delta_{C,M,\sigma^2,r}\ \lesssim \sigma^2rM\big(\mathrm{I}+\mathrm{II}\big).
\end{equation}
{\sc (The case $i=2$)} We assume subsequently that $\lambda_r>\lambda_{r+1}$, because otherwise $\mathrm{III}=\infty$ and the result follows with (\ref{eq: general1}).  We proceed similar to the case $i=1$ above, but with $k_{0,2}:=\infty$ and the auxiliary integer
\begin{equation*}
k_2^*\ :=\  \arg\max_{k\in\N}\bigg\{(\lambda_r^2-\lambda_{r+1}^2)\tilde{\Delta}_k-\Delta_k\geq\frac{1}{2}(\lambda_r^2-\lambda_{r+1}^2)\tilde{\Delta}_k\bigg\},
\end{equation*}
where $k_2^*:=0$ if this relation never holds. The sum
$$
\sum_{k < k_{0,2}}\E\Big(0\vee \big(W_{k,2}^0-\Omega_{k,2}\big)\Big)\ \leq\ \sum_{k\in\N}\E\Big(0\vee \big(W_{k,2}^0-(\lambda_r^2-\lambda_{r+1}^2)\tilde{\Delta}_k\big)\Big)
$$
can be treated analogously to the case $i=1$, with $\lambda_r^2-\lambda_{r+1}^2$ in place of $\lambda_r^2$. Similarly to (\ref{eq: k*}), $k>k_2^*$ implies
$$
c\sigma\lambda_1\tilde{\Delta}_k^{1/2}\sqrt{r(M-r)}>\tilde{\Delta}_k(\lambda_r^2-\lambda_{r+1}^2)/2.
$$
Since $\E(0\vee W_{k,2}^0)\leq \Delta_k$, the representation for the tail of the geometric series yields,  as above,
\begin{align*}
\sum_{k>k_2^*}\E \big(0\vee W_{k,2}^0\big)\ &\lesssim\ \sum_{k>k_2^*}\tilde{\Delta}_k^{1/2}\sigma\lambda_1\sqrt{r(M-r)}\ \lesssim \ \sigma^2rM\frac{\lambda_1^2}{{\lambda}_r^2-\lambda_{r+1}^2}.
\end{align*}
Combining case (i) and (ii) yields the proof of the Theorem.\hfill$\square$

\section{Proof of Theorem \ref{thm: lower bounds}}\label{sec: lower bounds 3}
Recall the definition $\tilde{\Delta}_k:= 2^{-k+2}r$, $k\in\N$. Since $\Arrowvert \pi_s(\id -\tilde{\pi}_s)\Arrowvert_{S_2}^2=\Arrowvert (\id -\pi_s)\tilde{\pi}_s\Arrowvert_{S_2}^2$ and $\pi_s-\tilde{\pi}_s=(\id-\pi_s)\tilde{\pi}_s-\pi_s(\id -\tilde{\pi}_s)$ is an orthogonal decomposition, observe first that 
\begin{align*}
\Arrowvert\pi_s C_{\alpha,s}\Arrowvert_{S_2}^2-\Arrowvert \tilde{\pi}_sC_{\alpha,s}\Arrowvert_{S_2}^2\ &=\ \alpha^2\Big(\Arrowvert\pi_s\Arrowvert_{S_2}^2-\Arrowvert \tilde{\pi}_s\pi_s\Arrowvert_{S_2}^2\Big)\\
&=\ \alpha^2\Arrowvert (\id -\tilde{\pi}_s)\pi_s\Arrowvert_{S_2}^2\\
&=\  \alpha^2\Arrowvert \tilde{\pi}_s(\id-\pi_s)\Arrowvert_{S_2}^2 =\ \frac{\alpha^2}{2}\Arrowvert \tilde{\pi}_s-\pi_s\Arrowvert_{S_2}^2,
\end{align*}
that is, 
$$\E\Arrowvert \pi_s(C_{\alpha,s}+\mathrm{E})\Arrowvert_{S_2}^2-\E\Arrowvert \tilde{\pi}_s(C_{\alpha,s}+\mathrm{E})\Arrowvert_{S_2}^2\in (\alpha^2\tilde{\Delta}_{k+1},\alpha^2\tilde{\Delta}_k] \ \Leftrightarrow\  \Arrowvert \tilde{\pi}_s-\pi_s\Arrowvert_{S_2}^2\in(2\tilde{\Delta}_{k+1},2\tilde{\Delta}_k].
$$
Recall the definition $$\GG_{M,s}(\delta, C)\ =\ \Big\{\tilde{\pi}_s\in\SS_{M,s}: \Arrowvert \pi_s C\Arrowvert_{S_2}^2-\Arrowvert \tilde{\pi}_sC\Arrowvert_{S_2}^2\leq\delta\Big\}.
$$
Note at this point that with $\bar{C}_{\alpha,M-s}:=\alpha U(\id-\id_s)V'$, we have  $\tilde{\pi}_r\in\GG_{M,r}(\delta, C_{\alpha,r}) \Leftrightarrow  (\id-\tilde{\pi}_r)\in\GG_{M,M-r}(\delta, \bar{C}_{\alpha,M-r}).
$ 
Define 
$$
 k^{**}\ :=\ \underset{k\in\N}{\arg\max}\Big\{{\alpha}^2 \tilde{\Delta}_{k}\geq d\sigma^2s(M-s) \Big\}
$$
for some $d>0$ to be specified later, 
and let $$
\DD_s(\alpha)\ :=\ \GG_{M,s}\big(\alpha^2\tilde{\Delta}_{k^{**}+1},C_{\alpha,s}\big)\ \ \ \text{as well as}\ \ \ \bar{\DD}_s({\alpha})\ :=\ \GG_{M,s}\big(\alpha^2\tilde{\Delta}_{k^{**}+1},\bar{C}_{\alpha,s}\big).
$$ 
Note that $k^{**}\rightarrow\infty$ as ${\alpha}\rightarrow\infty$. 
 It holds that
\begin{align*}
\E&\left(\sup_{\tilde{\pi}_s\in\SS_{M,s}}\big\Arrowvert\tilde{\pi}_s(C_{\alpha,s}+\mathrm{E})\big\Arrowvert_{S_2}^2-\big\Arrowvert{\pi}_s(C_{\alpha,s}+\mathrm{E})\big\Arrowvert_{S_2}^2\right)\\ 
&\ \ \ \ \ \ \ \ \ \geq\ \E\left(\sup_{\tilde{\pi}_s\in\DD_s(\alpha)}\big\Arrowvert\tilde{\pi}_s\mathrm{E}\big\Arrowvert_{S_2}^2-\big\Arrowvert{\pi}_s\mathrm{E}\big\Arrowvert_{S_2}^2 +2\alpha\,\tr\Big(\mathrm{E}'(\tilde{\pi}_s-\pi_s)\pi_s\Big)-ds(M-s)\sigma^2\right)\\
&\ \ \ \ \ \ \ \ \ \geq\ \E\left(\sup_{\tilde{\pi}_s\in\DD_s(\alpha)}2\alpha\,\tr\Big(\mathrm{E}'(\tilde{\pi}_s-\pi_s)\pi_s\Big)-ds(M-s)\sigma^2\right)\\ 
&\ \ \ \ \ \ \ \ \ \ \ \ \ \ \ \ \ \ -\  \E\left(\sup_{\tilde{\pi}_s\in\DD_s(\alpha)}\big\Arrowvert{\pi}_s\mathrm{E}\big\Arrowvert_{S_2}^2-\big\Arrowvert\tilde{\pi}_s\mathrm{E}\big\Arrowvert_{S_2}^2\right).
\end{align*}
Because of 
$\displaystyle
\underset{{\alpha}\rightarrow\infty}{\lim\sup}\,\E\left(\sup_{\tilde{\pi}_s\in\DD_s(\alpha)}\big\Arrowvert{\pi}_s\mathrm{E}\big\Arrowvert_{S_2}^2-\big\Arrowvert\tilde{\pi}_s\mathrm{E}\big\Arrowvert_{S_2}^2\right) = 0,
$
it remains to prove that
\begin{align}
\underset{{\alpha}\rightarrow\infty}{\lim\inf}&\,\max_{s\in\{r,M-r\}}\E\left(\sup_{\tilde{\pi}_s\in\DD_s(\alpha)}2\alpha\,\tr\Big(\mathrm{E}'(\tilde{\pi}_s-\pi_s)\pi_s\Big)-ds(M-s)\sigma^2\right)\nonumber\\
&=\ \sigma^2\,\underset{{\alpha}\rightarrow\infty}{\lim\inf}\,\max_{s\in\{r,M-r\}}\E\left(\sup_{\tilde{\pi}_s\in\DD_s(\alpha)}2({\alpha}/\sigma)\,\tr\big((\mathrm{E}/\sigma)'\tilde{\pi}_s\pi_s\big)-ds(M-s)\right)\label{eq: lower bound}\\
& \gtrsim\ \sigma^2r(M-r).\nonumber
\end{align}
First, we have
\begin{align*}
\E&\sup_{\tilde{\pi}_s\in \DD_s(\alpha)}\tr\big((\mathrm{E}/\sigma)'\tilde{\pi}_s\pi_s\big)\\
&=\ \E\sup_{\tilde{\pi}_s\in \DD_s(\alpha)}\tr\big((\mathrm{E}/\sigma)'(\id-\tilde{\pi}_s)\pi_s\big)\\
&=\ \E \sup_{\tilde{\pi}_s\in \DD_s(\alpha)}\Big(\tr\big((\mathrm{E}/\sigma)'\tilde{\pi}_s\big)-\tr\big((\mathrm{E}/\sigma)'\tilde{\pi}_s(\id-\pi_s)\big)\Big),%\\
%&\geq\ \E\sup_{\tilde{\pi}_r\in \DD(\tilde{\alpha})}\tr\big((\mathrm{E}/\sigma)'\tilde{\pi}_r\big)\ -\ \E\sup_{\tilde{\pi}_r\in \DD(\tilde{\alpha})}\tr\big((\mathrm{E}/\sigma)'\tilde{\pi}_r(\id-\pi_r)\big).
\end{align*}
which implies
\begin{align}
\E\sup_{\tilde{\pi}_s\in \DD_s(\alpha)}\tr\big((\mathrm{E}/\sigma)'&(\id-\tilde{\pi}_s)\pi_s\big)\ +\ \E\sup_{\tilde{\pi}_s\in \DD_s(\alpha)}\tr\big((\mathrm{E}/\sigma)'\tilde{\pi}_s(\id-\pi_s)\big)\label{eq: half}\\
& \geq\ \E\sup_{\tilde{\pi}_s\in \DD_s(\alpha)}\tr\big((\mathrm{E}/\sigma)'\tilde{\pi}_s\big).\nonumber
\end{align}
In  case $M\in 2\N$ and $s=M/2$,  since $\Arrowvert\pi_s-\tilde{\pi}_s\Arrowvert_{S_2}^2=\Arrowvert(\id-\pi_s)-(\id-\tilde{\pi}_s)\Arrowvert_{S_2}^2$, both expectations on the LHS in the inequality (\ref{eq: half}) are identical for reasons of symmetry, which leads to 
$$
\E\sup_{\tilde{\pi}_s\in \DD_s(\alpha)}\tr\big((\mathrm{E}/\sigma)'\tilde{\pi}_s\pi_s\big)\ \geq\ \frac{1}{2}\E\sup_{\tilde{\pi}_s\in \DD_s(\alpha)}\tr\big((\mathrm{E}/\sigma)'\tilde{\pi}_s\big)\ \ \ \text{in case $s=M/2$.}
$$
Although the polar decomposition of $(\id-\tilde{\pi}_s)\pi_s$ and $\tilde{\pi}_s(\id-\pi_s)$, respectively, suggests a similar symmetry argument% - in particular it shows that the singular values of $(\id-\tilde{\pi}_s)\pi_s$ and $\tilde{\pi}_s(\id-\pi_s)$ are the same
, we do not have a rigorous treatment of an argument of this type yet, and remain therefore with the inequality
\begin{align*}
\max_{s\in\{r,M-r\}}\E\sup_{\tilde{\pi}_s\in \DD_s(\alpha)}\tr\big((\mathrm{E}/\sigma)'(\id-\tilde{\pi}_s)\pi_s\big)\ &\geq\ \frac{1}{2}\,\E\sup_{\tilde{\pi}_r\in \DD_r(\alpha)}\tr\big((\mathrm{E}/\sigma)'\tilde{\pi}_r\big)%\\
%&=\ \frac{1}{2}\,\E\sup_{\tilde{\pi}_{M-r}\in \bar{\DD}_{M-r}(\tilde{\alpha})}\tr\big((\mathrm{E}/\sigma)'\tilde{\pi}_{M-r}\big)\\
\end{align*}
only. Note at this point that  
\begin{align*}
\E\sup_{\tilde{\pi}_r\in \DD_r(\alpha)}\tr\big((\mathrm{E}/\sigma)'\tilde{\pi}_r\big)\ &=\ \E\sup_{\tilde{\pi}_{M-r}\in \bar{\DD}_{M-r}({\alpha})}\tr\big((\mathrm{E}/\sigma)'\tilde{\pi}_{M-r}\big)\\ 
&=\ \E\sup_{\tilde{\pi}_{M-r}\in {\DD}_{M-r}({\alpha})}\tr\big((\mathrm{E}/\sigma)'\tilde{\pi}_{M-r}\big),
\end{align*}
where the second equality follows from invariance of the above expression under orthogonal transformation. By Sudakov's minoration and the  bound (\ref{eq: S2-2}) of Lemma \ref{lemma: bounds},
\begin{align*}
\E\left(\sup_{\tilde{\pi}_s\in\DD_s(\alpha)}\tr\big((\mathrm{E}/\sigma)'\tilde{\pi}_s\big)\right)\ &\gtrsim\ \delta\sqrt{\log N\big(\DD_r(\alpha), d_{S_2}, \delta\big)}\ \geq\ \delta\sqrt{r(M-r)}\sqrt{\log\Big(\frac{c\tilde{\Delta}_{k^{**}}^{1/2}}{\sqrt{2}\delta}\Big)}
\end{align*}
for any arbitrary $0<\delta<c\tilde{\Delta}_{k^{**}}^{1/2}/\sqrt{2}$, $s\in\{r,M-r\}$, where we used that $$
\E\left(\sup_{\tilde{\pi}_r\in\DD_r(\alpha)}\tr\big((\mathrm{E}/\sigma)'(\id-\tilde{\pi}_r)\big)\right)\ =\ \E\left(\sup_{\tilde{\pi}_r\in\DD_r(\alpha)}\tr\big((\mathrm{E}/\sigma)'\tilde{\pi}_r\big)\right)$$
and
$$
N\big(\DD_r(\alpha), d_{S_2}, \delta\big)\ \geq\ \bigg(\frac{c\tilde{\Delta}_{k^{**}}^{1/2}}{\sqrt{2}\delta}\bigg)^{r(M-r)}
$$
with the constant $c$ of Lemma \ref{lemma: bounds}. The choice 
 $\delta = c\tilde{\Delta}_{k^{**}}^{1/2}/8$ yields finally
\begin{equation}\label{eq: last}
 \E\left(\sup_{\tilde{\pi}_r\in\DD_r(\alpha)}\tr\big((\mathrm{E}/\sigma)'\tilde{\pi}_r\big)\right)\ \geq\  K\tilde{\Delta}_{k^{**}}^{1/2}\sqrt{r(M-r)}
\end{equation}
for some constant $K>0$ which does not depend on $\sigma $, $M$ and $\alpha$. 
Thus,
\begin{align}
(\ref{eq: lower bound})\ &=\ \sigma^2\,\underset{{\alpha}\rightarrow\infty}{\lim\inf}\max_{s\in\{r,M-r\}}\,\E\left(\sup_{\tilde{\pi}_s\in\DD_s(\alpha)}2({\alpha}/\sigma)\,\tr\big((\mathrm{E}/\sigma)'\tilde{\pi}_s\pi_s\big)-
dr(M-r)\right)\nonumber\\ 
&\geq\ \sigma^2\,\underset{{\alpha}\rightarrow\infty}{\lim\inf}\,\E\left(\sup_{\tilde{\pi}_r\in\DD_r(\alpha)}({\alpha}/\sigma)\,\tr\big((\mathrm{E}/\sigma)'\tilde{\pi}_r\big)-dr(M-r)\right).\label{eq: last 2}
\end{align}  
Choosing now $d$ in the definition of $k^{**}$ largest possible such that $
K\sqrt{d}-d \geq  K\sqrt{d}/2$ 
and plugging the lower bound (\ref{eq: last}) into (\ref{eq: last 2}) proves the Theorem.\hfill $\square$

\bibliography{donsker_bib}
\bibliographystyle{apalike}

\end{document}